\documentclass[ hidelinks, onefignum, onetabnum]{siamart220329}

\usepackage{amsmath, amssymb}
\usepackage{pgfplots}
\pgfplotsset{compat=1.18}
\usetikzlibrary{patterns}
\usepackage{lipsum}
\usepackage{graphicx}
\graphicspath{ {./images/} }
\usepackage{epstopdf}
\usepackage{mathtools}
\usepackage{multirow}
\usepackage{cancel}
\usepackage{ulem}
\usepackage{algorithmic}
\ifpdf
  \DeclareGraphicsExtensions{.eps,.pdf,.png,.jpg}
\else
  \DeclareGraphicsExtensions{.eps}
\fi

\makeatletter
\newcommand*{\addFileDependency}[1]{
  \typeout{(#1)}
  \@addtofilelist{#1}
  \IfFileExists{#1}{}{\typeout{No file #1.}}
}
\makeatother

\newsiamremark{remark}{Remark}
\newsiamremark{hypothesis}{Hypothesis}
\crefname{hypothesis}{Hypothesis}{Hypotheses}
\newsiamthm{claim}{Claim}

\newcommand{\N}{\mathbb{N}^{+}}
\newcommand{\Z}{\mathbb{Z}}
\newcommand{\R}{\mathbb{R}}
\newcommand{\C}{\mathbb{C}}

\newcommand{\imag}{{\mathrm{i}}}
\newcommand{\e}{{\mathrm{e}}}
\DeclareMathOperator*{\argmin}{argmin}
\newcommand{\A}[1]{\argmin_{#1}}
\DeclareMathOperator*{\expectation}{\mathbb{E}}
\newcommand{\E}[1]{\expectation_{#1}}

\newcommand{\hilb}{\mathfrak{H}}

\newcommand{\order}[1]{\mathcal{O}(#1)}

\newcommand{\totTime}{T}
\newcommand{\timeVar}{t}
\newcommand{\timeDisc}{N}
\newcommand{\timeInd}{n}

\newcommand{\stepSize}{h}

\newcommand{\spaVar}{x}
\newcommand{\spaDisc}{M}
\newcommand{\spaInd}{m}
\newcommand{\stateVar}{u}
\newcommand{\difEq}{f}
\newcommand{\trueState}[1]{\stateVar({#1})}
\newcommand{\hamil}{H}

\newcommand{\pot}{\widehat{V}}
\newcommand{\potOpr}{V}
\newcommand{\lapOpr}{\Delta}
\newcommand{\lap}{\widehat{\lapOpr}}
\newcommand{\lapdiag}{\lap_{\mathrm{diag}}}
\newcommand{\exactFlow}[1]{\psi_{#1}}

\newcommand{\numState}[1]{\stateVar_{#1}}
\newcommand{\numFlow}[1]{\Psi_{#1}}
\newcommand{\numFlowStep}[2]{\Psi_{#1, #2}}
\newcommand{\allParams}{\gamma}

\newcommand{\errorConst}{C}
\newcommand{\errorOrder}{p}

\newcommand{\difEqComp}[1]{\difEq^{[{#1}]}}
\newcommand{\exactFlowComp}[2]{\exactFlow{#2}^{[{#1}]}}
\newcommand{\numFlowComp}[3]{\numFlowStep{#2}{#3}^{[{#1}]}}
\newcommand{\numFlowCompOneStep}[2]{\numFlow{#2}^{[{#1}]}}
\newcommand{\potParams}{\alpha}
\newcommand{\kinParams}{\beta}
\newcommand{\splitDisc}{K}
\newcommand{\splitInd}{k}

\newcommand{\stateVarDist}{\mathcal{\MakeUppercase{\stateVar}}}

\newcommand{\paramTrans}{g}
\newcommand{\mlLossFn}{\mathcal{L}}
\newcommand{\mlParams}{\theta}
\newcommand{\mlFunc}[2]{\Phi(#1; #2)}

\newcommand{\cost}{\mathcal{C}}
\newcommand{\costPot}{\mathcal{C}^{[1]}}
\newcommand{\costKin}{\mathcal{C}^{[2]}}
\newcommand{\costConst}{C'}

\headers{Learning efficient and provably convergent splitting methods}{L. M. Kreusser, H. E. Lockyer, E. H. Müller, and P. Singh}

\title{Learning efficient and provably convergent splitting methods\thanks{Submitted to the editors \today.
\funding{Henry Lockyer is supported by a scholarship from the EPSRC Centre for Doctoral Training in Statistical Applied Mathematics at Bath (SAMBa), under the project EP/S022945/1.}}}

\author{Lisa M. Kreusser\footnotemark[3] \thanks{(Corresponding author).} 
\and Henry E. Lockyer\thanks{Department of Mathematical Sciences, University of Bath, Bath, BA2 7AY, United Kingdom
  (\email{lmk54@bath.ac.uk}, \email{hl785@bath.ac.uk}, \email{em459@bath.ac.uk}, \email{ps2106@bath.ac.uk}).}
\and Eike H. Müller\footnotemark[3]
\and Pranav Singh\footnotemark[3]}

\usepackage{amsopn}

\ifpdf
\hypersetup{
  pdftitle={Learning efficient and provably convergent splitting methods},
  pdfauthor={L. M. Kreusser, H. E. Lockyer, E. H. Müller, and P. Singh},
  colorlinks=True,
  allcolors=red!50!black
}
\fi

\begin{document}
\maketitle

\begin{abstract}
Splitting methods are widely used for solving initial value problems (IVPs) due to their ability to simplify complicated evolutions into more manageable subproblems. These subproblems can be solved efficiently and accurately, leveraging properties like linearity, sparsity and reduced stiffness. Traditionally, these methods are derived using analytic and algebraic techniques from numerical analysis, including truncated Taylor series and their Lie algebraic analogue, the Baker--Campbell--Hausdorff formula. These tools enable the development of high-order numerical methods that provide exceptional accuracy for small timesteps. Moreover, these methods often (nearly) conserve important physical invariants, such as mass, unitarity, and energy.
However, in many practical applications the computational resources are limited. Thus, it is crucial to identify methods that achieve the best accuracy within a fixed computational budget, which might require taking relatively large timesteps. In this regime, high-order methods derived with traditional methods often exhibit large errors since they are only designed to be asymptotically optimal. Machine Learning techniques offer a potential solution since they can be trained to efficiently solve a given IVP with less computational resources. However, they are often purely data-driven, come with limited convergence guarantees in the small-timestep regime and do not necessarily conserve physical invariants.
In this work, we propose a framework for finding machine learned splitting methods that are computationally efficient for large timesteps and have provable convergence and conservation guarantees in the small-timestep limit.
We demonstrate numerically that the learned methods, which by construction converge quadratically in the timestep size, can be significantly more efficient than established methods for the Schr\"{o}dinger equation if the computational budget is limited.
\end{abstract}


\section{Introduction} \label{sect:intro}
In this paper we consider first order initial value problems (IVPs) which arise in a wide range of physical applications. They encompass systems of ordinary differential equations (ODEs) with a finite dimensional state space, as well as the more general case of partial differential equations (PDEs) expressed as ODEs on an infinite dimensional Hilbert space $\hilb$. Mathematically, first order IVPs can be written as
\begin{equation} \label{eq:diffEq}
    \dot{\stateVar}(\timeVar) = \difEq(\stateVar, \timeVar), \quad \trueState{0} = \numState{0}, \quad \trueState{\timeVar} \in \hilb, \quad \timeVar \in [0,\totTime],
\end{equation}
for some $\totTime>0$, where $\stateVar$ is some state that evolves in time $\timeVar$ from a given initial state $\numState{0}$ under the action of a vector field $\difEq$ and $\dot{\stateVar}$ denotes the derivative of $u$ with respect to time $\timeVar$.  In general, it is not possible to solve IVPs like \eqref{eq:diffEq} analytically; even where closed form solutions do exist, their evaluation is often prohibitively expensive computationally. Thus we require numerical methods that are stable, accurate and computationally efficient. These are typically realised in terms of time-stepping methods where, for the sake of simplicity, we consider the evolution to the final time $\totTime = \timeDisc \stepSize$ as being split into $\timeDisc$ equal steps of size $\stepSize \ll \totTime$. A one-step numerical method is then uniquely defined by the \textit{forward map} from the approximate solution at a given time $\timeVar$ to the approximate solution at time $\timeVar + \stepSize$.

Splitting and composition methods \cite{blanes2024a} allow the separation of the evolution in \eqref{eq:diffEq} into simpler IVPs which can be solved efficiently and accurately; as a consequence they have been used successfully in many areas, see e.g.\ \cite{omelyan2003a}.
In particular, we assume that the vector field $\difEq = \difEqComp{1} + \difEqComp{2}$ can be split into two components $\difEqComp{1}$ and $\difEqComp{2}$, each of which defines an IVP that is simpler to solve numerically; as shown in \cite{glowinski2017a} this is indeed the case for a wide range of applications. 

Traditionally, splitting and composition methods methods are derived using analytic and algebraic conditions, including truncated Taylor series and their Lie algebraic analogue, the Baker–Campbell–Hausdorff (BCH) formula \cite{iserles2001a, bonfiglioli2012a} to guarantee consistency (or local error) of high order. Once stability is ensured, this leads to high-order convergence (of the global error) in the asymptotic limit of small timesteps $\stepSize \rightarrow 0$, see e.g.\ \cite{suzuki1990a, mcLachlan2002a, omelyan2003a, blanes2008a}. However, the asymptotic convergence of traditional methods implies that they are most efficient for small timestep $\stepSize$, which correspond to significant computational cost. It is also not obvious that a method which is designed to be fast in the limit $\stepSize \rightarrow 0$ will be the best choice for larger $\stepSize$. To make these points explicit we define two criteria to characterise an optimal numerical method:
\begin{description} \label{crit:optim}
    \item[C1:] The method has the fastest decrease in error in the asymptotic regime $\stepSize \to 0$.
    \item[C2:] The method has the smallest error for a given limited computational budget.
\end{description}
Traditionally, the focus has been on constructing methods that satisfy criterion \textbf{C1}. This, however, is not helpful if computational resources are limited and simulations have to be carried out at relatively large timestep sizes. This scenario, which is of significant practical interest, is the case we consider in this work. Our aim is to construct numerical methods with machine learning techniques such that they satisfy criterion \textbf{C2}, while being provably convergent in the limit $\stepSize \rightarrow 0$.

The exact solutions of IVPs often conserve a range of quantities that are known as invariants or first integrals. For example, for the Schrödinger equation the norm of the complex-valued solution does not change with time, which physically translates to the conservation of total probability, see~\cite{sakurai1985a}. More generally, for the Schrödinger equation the evolution of the IVP preserves the inner product of the underlying Hilbert space, a property known as unitarity, as well as the total energy, see \cite{sakurai1985a}. Similarly, Hamiltonian systems conserve a symplectic two-norm in phase space, see \cite{hairer2006a}. Such invariants often have a physical interpretation and are typically related to conservation laws. It is highly desirable to construct numerical methods which preserve these invariants at least approximately. This usually also improves the stability of the method since it limits the trajectory to a sub-manifold of the Hilbert space. An important example is the class of symplectic integrators, see e.g.\ \cite{hairer2006a}.

Machine Learning (ML) methods offer a promising alternative to traditional numerical methods for determining a numerical solution to an IVP, see e.g.\ \cite{karniadakis2021a} for a recent review. 
The key idea is to learn a forward solution map by training on a large set of initial and final values for the IVP. Provided the examples used for training are sufficiently representative of the (distribution of) initial conditions we are interested in solving, these methods can achieve high accuracy on unseen initial conditions  from the same distribution and potentially also generalise to a wider class of IVPs.
In their simplest form, ML methods are purely data-driven. They learn a fixed-cost forward solution map which is typically represented by a neural network. Unless measures are taken to specifically enforce inductive priors (see e.g. \cite{raissi2019a, offen2021a, müller2023a}), ML techniques do not guarantee the conservation of invariants such as unitarity and symplecticity which can be enforced with traditional splitting methods. 

With few exceptions such as \cite{chen2018a}, ML methods suffer from their detachment from the differential equation framework: since they do not parametrise the forward map as a function of the timestep size $\stepSize$, there is no sense in which convergence in the limit $\stepSize \to 0$ can be quantified. In contrast to traditional numerical methods which can be made more accurate by reducing the timestep size, it is not possible to achieve higher accuracy for these ML methods in a controlled way.
Even ML methods that learn the forward map with a variable timestep size $\stepSize$ do not necessarily generalise to different timestep sizes. In particular they may fail to converge in the limit $\stepSize \rightarrow 0$. Moreover, the stability of the forward map -- a crucial component for establishing the convergence -- is exceptionally hard to guarantee. 

However, provided sufficient training data is available, the model is sufficiently expressive, and the loss function weights large time step performance sufficiently, ML based methods have the potential to result in smaller overall errors, and thus superior performance for specific computational budgets and larger timesteps (in the sense of criterion \textbf{C2}) than traditional methods which aim to be efficient in the asymptotic regime (in the sense of criterion \textbf{C1}), see \cite{karniadakis2021a, grossmann2023a}. 
\subsection{Contributions} \label{sect:contributions}
In this paper we combine techniques from numerical analysis and ML to find efficient splitting schemes in the sense of criterion \textbf{C2} while being provably convergent in the limit $\stepSize \rightarrow 0$. Specifically we use machine learning to find splitting schemes that are tailored for a specific distribution of initial conditions, yet maintain many of the advantages of classical numerical methods, such as interpretability, generalisability, convergence, and conservation. As our numerical results show, the learned methods also generalise to scenarios not contained in the training set.

More specifically, we use ML methods to learn coefficients of splitting methods, while algebraically enforcing conditions that guarantee desirable properties such as consistency, stability and reversibility. We demonstrate the utility of this framework by learning splittings of medium to long length which result in lower errors than classical methods for large step sizes $\stepSize$. Restricting our search to physically plausible methods allows us to search a significantly lower dimensional submanifold. This reduces the training time compared to naive black-box ML approaches while also guaranteeing convergence. We numerically verify that the learned methods are efficient for the Schr\"{o}dinger equation with a double well potential and also demonstrate that they have the desirable conservation properties. While by construction the learned methods are formally only quadratically convergent in $\stepSize$, we also show that we can learn near fourth order methods.
\section{Traditional numerical analysis methods} \label{sect:numAna}
We restrict our attention to well-posed IVPs, using the definition of well-posedness found in \cite{evans2010a}. Hence any IVP of the form \eqref{eq:diffEq} considered in this work has a unique solution $\trueState{\timeVar}$ for all $\timeVar \in [0,\totTime]$ and is differentiable with respect to time and initial conditions. Arguments about existence and uniqueness of strong solutions to PDEs of the form \eqref{eq:diffEq} are active areas of research and are beyond the scope of this paper. Specifically for the linear Schrödinger equation existence and uniqueness of solutions under time-varying potentials are discussed in \cite{yajima1987a, ruggenthaler2015a}, while for the case of time-independant potentials the theory of one-parameter semi-groups can be applied, see \cite{engel1999a}. Once these PDEs are semi-discretised into a system of ODEs with the method of lines, the Picard-Lindelöf theorem provides conditions on the existence of unique solutions, see for instance \cite{giordano2022a}. In the case of time-independent potentials a computable closed-form expression for this unique solution is provided by the matrix exponential. However, evaluating this directly can be prohibitively expensive in practice \cite{moler2003a}. 
\subsection{Flows}
As we have assumed the well-posedness of our IVPs, we know there exits a unique solution $\trueState{\timeVar}$ for all times $\timeVar\in[0,\totTime]$. For a given time difference $\stepSize$ and vector field $\difEq$ we define the analytic flow $\exactFlowComp{\difEq}{\stepSize}:\hilb \to \hilb$ as the solution map
\begin{equation} \label{eq:trueFlow}
    \exactFlowComp{\difEq}{\stepSize}(\trueState{\timeVar}) = \trueState{\timeVar + \stepSize},
\end{equation}
which maps a state $\trueState{\timeVar}$ to its unique evolution $\trueState{\timeVar + \stepSize}$ after time $\stepSize$ under \eqref{eq:diffEq}. As the analytic flow $\exactFlowComp{\difEq}{\stepSize}$ encodes the true solutions of a differential equation, its precise form is typically unavailable or computationally infeasible. 
Numerical time-stepping methods, see for example \cite{hairer2006a, hairer2008a}, seek approximations to the analytic flow by approximating the solution at discrete times,
\begin{equation} \label{eq:timeDisc}
    \numState{\timeInd} \approx \trueState{\timeVar_\timeInd}, \qquad \timeInd \in \{0,1,\ldots, \timeDisc\},
\end{equation}
where $\timeDisc \in \N$ is the number of timesteps of size $\stepSize = \frac{\totTime}{\timeDisc} \in \R^+$ and $\timeVar_\timeInd = \timeInd \stepSize$, with $\timeVar_0=0$ and $\timeVar_\timeDisc = \totTime$. We consider one-step methods where the numerical solution $\numState{\timeInd+1}$ is the result of applying a numerical method (which may involve multiple evaluations of $\difEq$, or its subcomponents $\difEqComp{1}$ and $ \difEqComp{2}$ in the case of splitting methods) to the previous state $\numState{\timeInd}$. This can be seen as a numerical \textit{forward map} $\numFlowCompOneStep{\difEq}{\stepSize}:\hilb \to \hilb$,
\begin{equation} \label{eq:oneStepNumFlow}
    \numFlowCompOneStep{\difEq}{\stepSize}(\numState{\timeInd}) = \numState{\timeInd+1}.
\end{equation}
Inspired by the definition of the analytic flow in equation~\eqref{eq:trueFlow} we define the numerical flow as the map that sends the numerical state $\numState{\nu}$ to its evolution under $\timeInd$ timesteps of the forward map in \eqref{eq:oneStepNumFlow},
\begin{equation}\label{eq:numFlow}
    \numFlowComp{\difEq}{\timeVar_\timeInd}{\stepSize}(\numState{\nu}) = \left(\numFlowCompOneStep{\difEq}{\stepSize}\right)^\timeInd (\numState{\nu}) = \numState{\nu + \timeInd}.
\end{equation}
Note that we use the upper case $\Psi$ for numerical flow and the lower case $\psi$ for analytical flow. Unlike the exact flow in \eqref{eq:trueFlow}, the numerical flow in \eqref{eq:numFlow} is only defined at the discrete time points $\timeVar_\timeInd$, for $\timeInd \in \{0,1,\ldots,\timeDisc\}$. A numerical method and its associated numerical flow $\numFlowComp{\difEq}{\timeVar_\timeInd}{\stepSize}$ have order of convergence $\errorOrder \in \N$ if for each sufficiently regular \cite{hairer2006a, blanes2024a} $\difEq$ and $\numState{0}$ there exists a constant $\errorConst(\difEq, \numState{0}, \totTime)$ independent of $\stepSize$ such that 
\begin{equation} \label{eq:order}
    \left\lVert \numFlowComp{\difEq}{\timeVar_\timeInd}{\stepSize}(\numState{0}) - \exactFlowComp{\difEq}{\timeVar_\timeInd}(\numState{0}) \right\rVert \leq \errorConst \stepSize^{\errorOrder},
\end{equation}
for all $\timeInd$ as $\stepSize \to 0$. In other words, the numerical method is convergent of order $\errorOrder$ and the global error is of the order $\order{\stepSize^{\errorOrder}}$. Convergence requires both consistency and stability which imply that
 the numerical and analytical flows agree up to the $\errorOrder - 1^\text{th}$ term in the Taylor expansion in $\stepSize$.
\subsection{Splittings}
For the differential equation of the form \eqref{eq:diffEq} we assume that 
the vector field, $\difEq = \difEqComp{1} + \difEqComp{2}$, can be split into two (sub-)components, $\difEqComp{1}$ and $\difEqComp{2}$. We further assume that the two IVPs, defined by the flow under the individual vector fields $\difEqComp{1}$ and $\difEqComp{2}$, are well-posed and have known analytic solutions which can readily be computed. The corresponding analytic flows are called the ``subflows'' and denoted by $\exactFlowComp{1}{\timeVar}$ and $\exactFlowComp{2}{\timeVar}$ respectively. A numerical method for the IVP in \eqref{eq:diffEq} can be constructed by interleaving the sub-flows. Two such well known splitting methods are
\begin{align}
    \text{Trotter} &: \numFlowCompOneStep{\difEq}{\stepSize} = \exactFlowComp{2}{\stepSize} \circ \exactFlowComp{1}{\stepSize} \qquad \text{and} \label{eq:trottSolve}
    \\
    \text{Strang} &: \numFlowCompOneStep{\difEq}{\stepSize} = \exactFlowComp{1}{0.5\stepSize} \circ \exactFlowComp{2}{\stepSize} \circ \exactFlowComp{1}{0.5\stepSize}. \label{eq:strangSolve}
\end{align}
Trotter and Strang are methods of order one and two respectively, in the sense of \eqref{eq:order}, see \cite{blanes2024a}. As a consequence, the Taylor coefficients of the analytical and numerical flows for all points in the temporal discretisations agree for the constant, $\stepSize^0$ term for both methods; for Strang the coefficients of the $\stepSize^1$ terms also agree. Both methods are part of the larger family of splitting schemes that can be parameterised as follows,
\begin{equation} \label{eq:splitParam}
    \numFlowComp{\difEq}{\totTime}{\stepSize}(\, \cdot \, ; \potParams, \kinParams) = \bigcirc_{\timeInd = 1}^{\timeDisc}\numFlowCompOneStep{\difEq}{\stepSize}(\, \cdot \, ; \potParams, \kinParams), \text{ where } \numFlowCompOneStep{\difEq}{\stepSize}(\, \cdot \, ; \potParams, \kinParams) = \bigcirc_{\splitInd = 1}^{\splitDisc} \exactFlowComp{2}{\kinParams_\splitInd \stepSize} \circ \exactFlowComp{1}{\potParams_\splitInd \stepSize}.
\end{equation}
Here, $\bigcirc$ represents repeated composition, i.e.\ $\bigcirc_{\timeInd = 1}^{\timeDisc}\varphi_\timeInd = \varphi_\timeDisc \circ \dots \circ \varphi_2 \circ \varphi_1$ and $\numFlowComp{\difEq}{\totTime}{\stepSize}(\, \cdot \,; \potParams, \kinParams)$, $\numFlowCompOneStep{\difEq}{\stepSize}(\, \cdot \,; \potParams, \kinParams)$,  $\exactFlowComp{1}{\potParams_\splitInd\stepSize}$ and $\exactFlowComp{2}{\kinParams_\splitInd\stepSize}$ map from $\hilb$ to $\hilb$. We call a splitting of the form \eqref{eq:splitParam} a $\splitDisc$-stage method with $\timeDisc$ individual steps and a timestep of size $\stepSize=\totTime/\timeDisc$. Observe that in general \eqref{eq:splitParam} requires a total of $2 \splitDisc \timeDisc$ sub-flow evaluations. The vector of $2\splitDisc$ parameters $ [\potParams, \kinParams] = [\potParams_1, \potParams_2, \dots, \potParams_\splitDisc, \kinParams_1, \kinParams_2, \dots, \kinParams_\splitDisc]$ continuously parameterises all\footnote{Technically in \eqref{eq:splitParam} we have imposed that the flow $\exactFlowComp{1}{\stepSize}$ associated with $\difEqComp{1}$ comes first, precluding some splittings, but this can be solved by relabelling $\difEqComp{1}\leftrightarrow\difEqComp{2}$.} possible splittings of maximum length $\splitDisc$.

We assume differentiability of the parameterised numerical flow $\numFlowComp{\difEq}{\totTime}{\stepSize}(\, \cdot \, ; \potParams, \kinParams)$ with respect to the parameters $\potParams_\splitInd$, $\kinParams_\splitInd$. To  ensure 
this, it suffices to assume the differentiability of the analytical sub-flows $\exactFlowComp{1}{\potParams_\splitInd \stepSize}$ and $\exactFlowComp{2}{\kinParams_\splitInd \stepSize}$ 
with respect to time and initial conditions which naturally follows from the assumed well-posedness of the sub-flows. 
\subsection{Order conditions}
While some particular IVPs allow exact splittings \cite{bernier2019a, bernier2020a}, in general splitting methods introduce errors which depend on the particular IVP that is solved. Nevertheless, some splitting schemes do have (near) universal desirable properties. For example, imposing so-called order conditions on the parameters $\potParams$ and $\kinParams$ ensures that the corresponding numerical methods have a particular order. The `first order consistency' order condition, henceforth simply called ``consistency'', is given by,
\begin{align} \label{eq:cons}
    \sum_{\splitInd = 1}^{\splitDisc} \potParams_\splitInd = 1 \ \text{and} \ \sum_{\splitInd = 1}^{\splitDisc}\kinParams_\splitInd = 1.
\end{align}
It can then be shown that \eqref{eq:cons} implies convergence with order $\errorOrder=1$ in the sense of \eqref{eq:order} if the scheme is also stable. Moreover, any  splitting scheme with convergence of order $\errorOrder\geq1$ necessarily needs to satisfy the consistency condition \eqref{eq:cons}. Thus, without loss of generalisation, we assume in the remainder of the manuscript that the consistency condition \eqref{eq:cons} is always satisfied.

\begin{figure}[tb]
    \centering
    \includegraphics[width=0.85\textwidth]{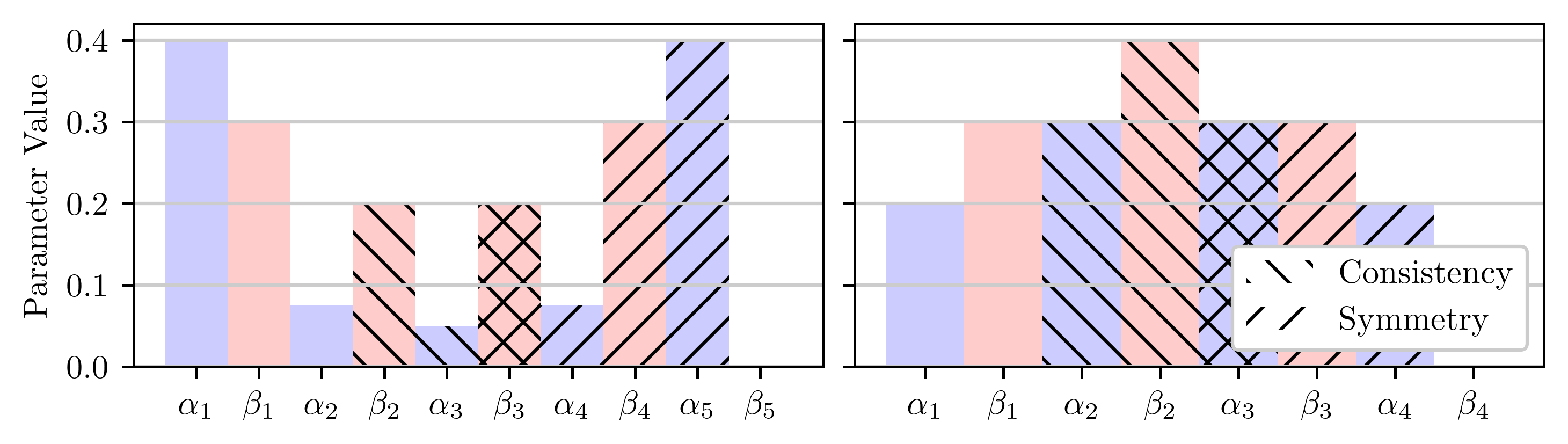}
    \caption{Visualisation of the $\potParams$ and $\kinParams$ for two symmetric and  consistent splitting schemes with $\splitDisc=5$ (left) and $\splitDisc=4$ (right) stages. Note how symmetry  constrains half the parameters, including the trailing $\kinParams_K$ that is symbolically zero. Consistency fixes two additional parameters to ensure that the $\potParams$ and $\kinParams$ sum to one as in \eqref{eq:cons}.
    \label{fig:paramTransform}}
\end{figure}

The simplest consistent splitting scheme with the smallest value of $\splitDisc$ is the Trotter method defined by the parameters $[\potParams_1, \kinParams_1] = [1.0, 1.0]$, which satisfy \eqref{eq:cons}. For longer splittings with larger $\splitDisc$ we can impose further order conditions on the parameters $\potParams_\splitInd$ and $\kinParams_\splitInd$. For example, enforcing symmetry under time-reversal ensures that the splitting error is an even function of the timestep size, i.e. all odd powers of $\stepSize$ in the Taylor expansion in $\stepSize$ vanish, and hence,
\begin{equation} \label{eq:symm}
    \numFlowComp{\difEq}{\totTime}{\stepSize} = \left(\numFlowComp{\difEq}{-\totTime}{-\stepSize}\right)^{-1} \ \Rightarrow \ \numFlowComp{\difEq}{\totTime}{\stepSize}(\stateVar) = \exactFlowComp{\difEq}{\totTime}(\stateVar) + \order{\stepSize^{\errorOrder}}, \quad \errorOrder \in 2\Z.
\end{equation}
When combined with stability, symmetry under time-reversal, henceforth simply called ``symmetry'', implies even order convergence. Since we assume consistency, symmetry ensures convergence of order at least $p=2$. As shown in \cite{blanes2024a}, symmetry can be guaranteed if the coefficient vectors that define the scheme in \eqref{eq:splitParam} are palindromic, i.e. of the form $(\potParams_1, \potParams_2, \dots, \potParams_\splitDisc) = (\potParams_\splitDisc, \potParams_{\splitDisc - 1}, \dots, \potParams_1)$ and 
$(\kinParams_1, \kinParams_2, \dots, \kinParams_{\splitDisc-1},0) = (\kinParams_{\splitDisc - 1}, \kinParams_{\splitDisc - 2}, \dots, \kinParams_1, 0)$ with the symbolic zero $\beta_K=0$.
Enforcing  consistency and symmetry imposes constraints on $\potParams$ and $\kinParams$ and hence reduces the number of degrees of freedom that define the splitting scheme: to parameterise all possible symmetric and consistent splittings we only require $\splitDisc - 2$ independent parameters rather than the $2\splitDisc$ parameters required to parameterise all splittings. Since these conditions are straightforward to impose and reduce the dimension of the search space, we restrict our attention to symmetric and consistent methods for the rest of the manuscript. The parameter of two such methods are visualised in Figure~\ref{fig:paramTransform}.
\subsection{Composing splitting methods} \label{sec:composing_splitting_methods}
The simplest, and indeed unique, splitting scheme of length two that is both consistent and symmetric (and hence second order) is Strang, defined by the parameters $[0.5, 0.5, 1.0, 0.0]$. However, finding higher-order conditions requires significant work. An alternative is to build up higher-order methods by composing multiple copies of lower order methods, see \cite[II.4]{hairer2006a}. An example of a composition method is the triple jump technique which we discuss in more detail in Appendix~\ref{sect:composition} and which can be applied to Strang to construct a splitting method of order four, referred to as Yoshida in this paper. A disadvantage with using composition methods to find new splittings is that they implicitly restrict the methods that can be found, and the number of stages grows exponentially with the desired order.
\subsection{Cost estimates}
The symbolic zero $\kinParams_\splitDisc=0$ in symmetric methods is worth remarking on as the costs for the evaluation of symmetric and general methods differs significantly. We  denote the cost of a single evaluation of the subflows $\exactFlowComp{1}{\potParams_\splitInd \stepSize}$, $\exactFlowComp{2}{\kinParams_\splitInd \stepSize}$ by $\costPot$ and $\costKin$ respectively.
For evaluating $\timeDisc$ individual steps of a general $\splitDisc$-stage method, the total cost is given by  $\timeDisc \splitDisc (\costPot+\costKin)$. For symmetric $\splitDisc$-stage methods, composing $\timeDisc$ individual steps does not require the evaluation of the identity flows associated with the symbolic zeros $\kinParams_\splitDisc=0$, while adjacent sub-flows associated with $\difEqComp{1}$ can be combined since $\exactFlowComp{1}{\potParams_1 \stepSize} \circ \exactFlowComp{1}{\potParams_\splitDisc \stepSize} = \exactFlowComp{1}{(\potParams_1 + \potParams_\splitDisc) \stepSize}$. For symmetric methods, the total cost therefore satisfies $\timeDisc (\splitDisc \costPot+(\splitDisc-1)\costKin)-(\timeDisc-1) \costPot = \timeDisc (\splitDisc-1) (\costPot+\costKin)+\costPot$, where the subtraction of $(\timeDisc-1) \costPot$ accounts for the reduction of the cost due to combination of sub-flows. Overall, this leads to relative cost savings of a factor $(\splitDisc-1)/\splitDisc+\mathcal{O}(\stepSize)$ per timestep. To summarise, the total cost of numerically integrating the IVP \eqref{eq:diffEq} to the final time $\totTime$ with a $\splitDisc$-stage splitting method and a timestep size of $\stepSize=\totTime/\timeDisc$ is given by,
\begin{equation}
    \cost = 
    \left. \begin{cases}
        \timeDisc  (\splitDisc-1)\left( \costPot + \costKin\right)  + \costPot & \text{if symmetric} \\
        \timeDisc  \splitDisc\left( \costPot + \costKin \right) & \text{in general}
    \end{cases} \right\}
    \ge \costConst(\splitDisc,\totTime) h^{-1}, 
    \label{eqn:total_cost_h}
\end{equation}
where the constant,
\begin{equation}
    \costConst(\splitDisc,\totTime)=
    \totTime\cdot
    \begin{cases}
        (\splitDisc-1)\left(\costPot + \costKin\right)  &\text{if symmetric,}\\
        \splitDisc\left( \costPot+\costKin\right) & \text{in general,}\\
    \end{cases}
\end{equation}
is independent of $\stepSize$ and the bound in \eqref{eqn:total_cost_h} is sharp up to corrections of relative $\mathcal{O}(h)$.
As a consequence, up to corrections of $\mathcal{O}(\stepSize)$, which are amortised for $\timeDisc\gg1$, a single step of symmetric Strang ($\splitDisc=2$) costs the same as a step of non-symmetric Trotter ($\splitDisc=1$). 
\subsection{Computational efficiency}
To choose an optimal numerical time-stepping method, we aim to minimise the computational cost $\cost=\cost(\varepsilon)$ for a given tolerance $\varepsilon$ on the numerical error $\left\lVert \numFlowComp{\difEq}{\totTime}{\stepSize}(\numState{0}) - \exactFlowComp{\difEq}{\totTime}(\numState{0}) \right\rVert$. If the method $\numFlowComp{\difEq}{\totTime}{\stepSize}(\numState{0})$ is convergent of order $\errorOrder$, then according to \eqref{eq:order} the numerical error at the final time is smaller than the tolerance $\varepsilon$ provided that the timestep size $\stepSize$ is chosen such that
\begin{equation}
    \errorConst(\difEqComp{1}, \difEqComp{2}, \numState{0}, \totTime;[\potParams, \kinParams]) h^p \le \varepsilon.\label{eqn:error_bound}
\end{equation}
Here we have made the dependence of the constant $\errorConst$ on the timestepping method explicit through the  parameters $[\potParams, \kinParams]$.
Combining \eqref{eqn:total_cost_h} and \eqref{eqn:error_bound}, we see that the total cost as a function of the tolerance $\varepsilon$ is bounded from below by
\begin{equation}
    \cost(\varepsilon) \ge \costConst(\splitDisc,\totTime) 
    \errorConst(\difEqComp{1}, \difEqComp{2}, \numState{0}, \totTime;[\potParams, \kinParams])^{1/p} \varepsilon^{-1/p}.\label{eqn:cost_bound}
\end{equation}
Provided that the bounds in \eqref{eqn:total_cost_h} and \eqref{eqn:error_bound} are sharp, this implies that in the limit $\varepsilon,\stepSize\to 0$ high-order methods will be most cost effective. However, in practical applications choosing a numerical method based on this argument (which corresponds to criterion \textbf{C1}) might be premature for two reasons: firstly, the constants that appear in \eqref{eqn:cost_bound} also depend on the timestepping method and might be of a similar size as the factor $\varepsilon^{-1/p}$ for moderate tolerances $\varepsilon$. Consider, for example, the discussion in Appendix~\ref{sect:composition}: composing splitting methods with the triple-jump technique leads to a $\mathcal{O}(3^p)$ exponential growth in the number of stages, $\splitDisc$, and thus to a similar increase of $\costConst(\splitDisc,\totTime)$. As a consequence, for a fixed computational budget (criterion \textbf{C2}) and thus moderately large step size $\stepSize$, a smaller error might be obtained with a lower-order splitting method with a small error constant $\errorConst(\difEqComp{1}, \difEqComp{2}, \numState{0}, \totTime;[\potParams, \kinParams])$ than by resorting to a higher-order method which may have a larger error constant. This effect is particularly pronounced when low to moderate levels of accuracy $\varepsilon$ are required or acceptable, a typical scenario in applications which are either significantly constrained by computational cost or where the IVP only provides a moderately accurate model of the underlying physical system. Hence, at any given order it is highly desirable to minimise the constant $\errorConst(\difEqComp{1}, \difEqComp{2}, \numState{0}, \totTime;[\potParams, \kinParams])$. 
\subsection{Identifying efficient splitting methods}
To a certain extent studies such as \cite{omelyan2003a} attempt to obtain a lower error constant by considering a large number of splittings, typically obtained by solving the algebraic order conditions, and identifying the splittings that feature the smallest coefficients accompanying the leading commutators terms in the residual BCH expansion. By design, this strategy is only asymptotically relevant, i.e. for small step size $\stepSize$. Moreover, it ignores the magnitude of the commutators and hence does not minimize the constant $\errorConst(\difEqComp{1}, \difEqComp{2}, \numState{0}, \totTime;[\potParams, \kinParams])$ that also depends on the particular IVP and its splitting given by $\difEqComp{1}, \difEqComp{2}$, the initial condition $\numState{0}$ and the final time $\totTime$. These factors are not straightforward to incorporate in an algebraically motivated approach.

When the subflows are unitary or symplectic, a splitting method guarantees the exact conservation of a shadow Hamiltonian or near conservation of the true Hamiltonian, and the error in energy is sometimes more relevant than the $L_2$-error of the solution. In the context of the Hamiltonian Monte Carlo sampling method, for instance, it becomes possible to obtain optimal splittings for larger step sizes by minimizing the expected error in energy under certain analytic priors on the problem \cite{sanzSerna2017a,blanes2014a}. 
\subsection{A data-driven approach for minimising the error constant}
Our goal is to construct a method which is robust in the sense that it produces a small error for a range of different vectors fields $\difEq$ and splittings $\difEqComp{1}, \difEqComp{2}$, as well as different initial conditions $\numState{0}$ and final times $\totTime$. For example, keeping $\difEqComp{1}, \difEqComp{2}$ and $\totTime$ fixed, we might be interested in,
\begin{equation}
    \overline{\errorConst}(\difEqComp{1}, \difEqComp{2}, \totTime;[\potParams, \kinParams]) = \sup_{\numState{0} \in \mathbb{U}}\left\{ \errorConst(\difEqComp{1}, \difEqComp{2}, \numState{0}, \totTime;[\potParams, \kinParams])\right\},
\end{equation}
 for all possible initial conditions $\numState{0}$   from some given set $\mathbb{U}$. This would allow replacing the specific bound in \eqref{eqn:error_bound} by the more general bound,
\begin{equation}
    \overline{\errorConst}(\difEqComp{1}, \difEqComp{2}, \totTime;[\potParams, \kinParams]) h^p \le \varepsilon\label{eqn:error_bound_general},
\end{equation} 
which holds for all $\numState{0} \in \mathbb{U}$. Similarly, we might want to make the bound problem-independent by considering the supremum over a range of $\difEq$ and splittings $\difEqComp{1}, \difEqComp{2}$. However, it is usually not possible to write down the dependence of $\errorConst$ (let alone $\overline{\errorConst}$) on $\difEqComp{1}, \difEqComp{2}, \numState{0}$ and $\totTime$ in closed form, so providing analytical bounds of the form in \eqref{eqn:error_bound_general} is typically not feasible. Instead, we could choose to quantify the error as the \textit{expected} error when drawing the initial condition $\numState{0}$ from some probability distribution $\stateVarDist$. In practice, we would then estimate the expected error of the distribution $\stateVarDist$ with a large, but finite sample $\mathbb{U}$ drawn from $\stateVarDist$. This is the approach which we will take in this work. As will be discussed in the following section, it fits naturally with the machine learning framework that we propose.
\section{Learned splitting schemes}
We now explain how we use Machine Learning (ML) to construct efficient splitting schemes which are guaranteed to have desirable properties such as consistency and symmetry. To do this, we first review some key concepts in ML in the context of the problem we consider in this paper. After this, we discuss how to enforce consistency and time reversal symmetry  through a parameter transformation and outline our training algorithm.
\subsection{ML notation and concepts} \label{sec:ML_notation_and_concepts}
In the setting considered here, a general supervised machine learning approach requires three ingredients:
\begin{itemize}
    \item A dataset $\mathbb{U}$, which consists of initial conditions for the IVP in \eqref{eq:diffEq}. For each $\numState{0} \in \mathbb U$, we assume that we can approximate the solution $\trueState{\totTime}=\exactFlowComp{\difEq}{\totTime}(\numState{0})$ of the IVP at time $\totTime$ to high precision, i.e.\ $\numState{\text{Ref}}\approx \trueState{\totTime}$, and we therefore construct labelled pairs $(\numState{0},\numState{\text{Ref}})$.
    \item A numerical flow function $\mlFunc{\, \cdot \, }{\mlParams}$ which encodes a numerical method and is parametrised by learnable parameters $\mlParams\in \mathbb R^d$ for some given $d\in \mathbb N$.
    \item A loss function $\mlLossFn(\mlParams)$ which measures how well $\mlFunc{\, \cdot \, }{\mlParams}$ approximates the analytic flow $\exactFlowComp{\difEq}{\totTime}$ for initial conditions $\numState{0} \in \mathbb{U}$.
\end{itemize}
During training, the parameters $\mlParams$ are tuned such that $\mlFunc{\, \cdot \, }{\mlParams}$ approximates $\exactFlowComp{\difEq}{\totTime}$. Ideally, we would like to choose $\mathbb U$ to be the set of \textit{all} possible physical initial conditions. However, this is typically impractical and in the case of PDEs, where this set is parametrised by an infinite number of degrees of freedom, impossible. Instead, we assume here that the initial conditions $\numState{0} \sim \stateVarDist$ are drawn from a probability distribution which we denote by $\stateVarDist$. We define the loss function $\mlLossFn(\mlParams)$  as the expectation of the mean squared error over the distribution $\stateVarDist$ which is given by,
\begin{equation} \label{eq:simpleMLloss}
    \mlLossFn(\mlParams) = \E{\numState{0} \sim \stateVarDist} \left\lVert \mlFunc{\numState{0}}{\mlParams} - \numState{\text{Ref}} \right\rVert_2^2.
\end{equation}
In practice, only a finite number of samples is considered so that the expectation in the loss function \eqref{eq:simpleMLloss} can be replaced by a sum. The optimal parameters $\mlParams^*$ are then obtained by minimising the loss, i.e. 
\begin{equation} \label{eq:simpleML}
    \mlParams^* = \A{\mlParams} \mlLossFn(\mlParams).
\end{equation}
Our aim is therefore to determine $\mlParams^*$ results in the learned function $\mlFunc{\, \cdot \,}{\mlParams^*}$ whose evaluation is efficient in the sense of \textbf{C2}, and which generalises to unseen inputs $\numState{0}$. This includes initial conditions with low probability or even $\numState{0}$ which are not in the distribution $\stateVarDist$. Usually this is possible if the function $\mlFunc{\, \cdot \,}{\mlParams}$ is both sufficiently expressive and enough samples $\numState{0}\sim \stateVarDist$ are used for training to find the optimal $\mlParams^*$ while avoiding overfitting. Overfitting can be mitigated by adding regularisers or by including inductive biases in the construction of the function $\mlFunc{\, \cdot \,}{\mlParams}$.

To determine $\mlParams^*$, brute force grid search or parameter sweep methods can be used to find an approximate minimiser of the loss function \eqref{eq:simpleMLloss} provided the parameters are contained in a bounded subset of $\mathbb{R}^d$ for some small $d$. For  large  $d$, these approaches are too expensive since the cost grows exponentially in $d$. For high-dimensional search spaces, $\mlParams^*$ is usually determined using gradient-based methods which only require differentiability of $\mlLossFn$ with respect to $\mlParams$. If $\numState{0} \sim \stateVarDist$ is randomly drawn for each gradient computation we call this stochastic optimisation (SO). Examples of SO include stochastic gradient descent (SGD) or improved variants like Adam \cite{kingma2014a}.
\subsection{Incorporating flows into ML}\label{sec:mlflows}
We let our learnable numerical flows $\mlFunc{\, \cdot \,}{\mlParams}$ be the splitting methods introduced in Section~\ref{sect:numAna}. Hence $\mlFunc{\, \cdot \,}{\mlParams}$ is obtained by composing multiple applications of the forward map $\numFlowCompOneStep{\difEq}{\stepSize}(\, \cdot \, ; \potParams, \kinParams)$ and the learnable parameters $\mlParams$ are the coefficients $\potParams, \kinParams$ that define the splitting scheme. In symbols, we define the learnable flow as,
\begin{equation} \label{eqn:learned_splitting_method}
    \mlFunc{\numState{0}}{\mlParams}=\numFlowComp{\difEq}{\totTime}{\stepSize}(\numState{0} ; \potParams, \kinParams):=\left(\numFlowCompOneStep{\difEq}{\stepSize}(\numState{0};\potParams, \kinParams)\right)^\timeDisc.
\end{equation}
By restricting the values of the learnable parameters $\potParams, \kinParams$ as in \eqref{eq:cons} we can guarantee the consistency order condition. Provided the IVP is well-defined and the method is stable, this will imply that the learned flow will produce provably convergent solutions in the limit $\stepSize\rightarrow 0$. Higher-order consistency can be enforced by imposing further order conditions, such as symmetry under time reversal. Other constraints from classical numerical analysis such as unitarity or symplecticity can also be enforced in this ML framework by enforcing the very same classical conditions from pre-existing numerical analysis. With the learned flow $\mlFunc{\numState{0}}{\mlParams}$ in \eqref{eqn:learned_splitting_method} the loss function in \eqref{eq:simpleMLloss} becomes,
\begin{equation} \label{eq:lossFn}
    \mlLossFn(\potParams, \kinParams) = \E{\numState{0} \sim \stateVarDist} \left\lVert \numFlowComp{\difEq}{\totTime}{\stepSize}(\numState{0} ; \potParams, \kinParams) - \numState{\text{Ref}}
    \right\rVert_2^2.
\end{equation}
By construction, the loss function in \eqref{eq:lossFn} is differentiable in $\potParams$ and $\kinParams$, as the numerical flow was differentiable in $\potParams$ and $\kinParams$. As a consequence, it can be minimised with gradient-based methods as discussed in Section~\ref{sec:ML_notation_and_concepts}. In analogy to \eqref{eq:simpleML}, we obtain the optimal parameters $\potParams^*, \kinParams^*$ as
\begin{equation} \label{eq:lossArgMin}
    \potParams^*, \kinParams^* = \A{\potParams, \kinParams} \mlLossFn(\potParams, \kinParams) \text{ for $\mlLossFn$ defined in~\eqref{eq:lossFn}}.
\end{equation}
The minima of the loss function \eqref{eq:lossArgMin} depends on the distribution $\stateVarDist$. In other words, the learned splittings are tailored to the specific distribution of initial values, and hence a specific distribution of IVPs. Minimizing the loss function~\eqref{eq:lossFn} minimises the expected error. For an order $\errorOrder$ splitting the error is of the form 
\begin{align}\label{eqn:error_expansion}
\left\lVert \numFlowComp{\difEq}{\totTime}{\stepSize}(\numState{0}; \potParams, \kinParams) - \numState{\text{Ref}}
\right\rVert &= \errorConst_0(\difEqComp{1}, \difEqComp{2}, \numState{0}, \totTime;[\potParams, \kinParams])\stepSize^{\errorOrder} \\
\nonumber & \qquad + \sum_{m=1}^\infty \errorConst_m(\difEqComp{1}, \difEqComp{2}, \numState{0}, \totTime;[\potParams, \kinParams])\stepSize^{\errorOrder+m}
\end{align}
Observe that for small timestep sizes $\stepSize$ the dominant contribution to the error will come from the first term in \eqref{eqn:error_expansion}. Minimizing the loss function~\eqref{eq:lossFn} will minimise the expected value of all error constants $\errorConst_m(\difEqComp{1}, \difEqComp{2}, \numState{0}, \totTime;[\potParams, \kinParams])$ and in particular the expected value of the leading constant $\errorConst_0(\difEqComp{1}, \difEqComp{2}, \numState{0}, \totTime;[\potParams, \kinParams])$. As a consequence, the error of the learned method is potentially significantly smaller than that of a method constructed with traditional techniques from numerical analysis, which rely purely on algebraic and analytic techniques, in particular for larger values of $\stepSize$ where asymptotic analysis typically breaks down. In other words, the learned method is efficient in the sense of criterion \textbf{C2}.
\subsection{Enforcing consistency and symmetry through parameter transformations}\label{sec:parameter_transform}
Training on the landscape defined by the loss function~\eqref{eq:lossFn} for a splitting method $\numFlowComp{\difEq}{\totTime}{\stepSize}(\numState{0} ; \potParams, \kinParams)$ is non-trivial since it is non-convex and ill-conditioned, in particular for small values of $\stepSize$. 
To see this, consider the expansion of the error in \eqref{eqn:error_expansion}. For sufficiently long splitting methods the manifold of all parameters $\potParams, \kinParams$ will likely contain a sub-manifold of parameters for which the method is of order $\errorOrder+1$, i.e.\ the constant $\errorConst_0(\difEqComp{1}, \difEqComp{2}, \numState{0}, \totTime;[\potParams, \kinParams])$ will be identically zero on this sub-manifold. Hence, varying the parameters within this sub-manifold will only lead to small changes in the loss function due to changes in the error constants $\errorConst_m(\difEqComp{1}, \difEqComp{2}, \numState{0}, \totTime;[\potParams, \kinParams])$ with $m\ge 1$. However, moving away from the sub-manifold will at least reduce the order of the method from $\errorOrder+1$ to $\errorOrder$, and thus incur changes of $\mathcal{O}(\stepSize^{-1})$ in the loss function. This argument can be applied recursively: sub-manifolds of higher-order methods are nested within manifolds that describe lower-order methods. As a consequence, the condition number of the Hessian  grows with some power of the inverse timestep size $\stepSize^{-1}$. Furthermore, this problem is likely to get worse for splitting methods with more sub-flows $\splitDisc$, i.e.\ higher-dimensional parameter spaces. To compensate for the poor conditioning, we would have to severely restrict the learning rate which renders SO methods very inefficient.

Given that the consistency order condition \eqref{eq:cons} is easy to enforce and required for provable convergence guarantees of the learned method, it would be unwise not to incorporate it. The consistency order condition \eqref{eq:cons} restricts the parameters $\potParams, \kinParams$ to an affine hyperplane.  Symmetric splitting methods  can be obtained by imposing an additional set of linear constraints on $\potParams, \kinParams$, which -- together with~\eqref{eq:cons} -- reduce the parameters to a lower-dimensional hyperplane that contains methods of even order.

To restrict ourselves to the  sub-manifold that enforces consistency and symmetry, we express the original $\potParams, \kinParams\in \mathbb R^\splitDisc$ in terms of suitable coordinates $\allParams_\potParams\in \mathbb R^{ \lfloor \frac{ \splitDisc  -1}{2}\rfloor}$, $\allParams_\kinParams \in \mathbb R^{ \lfloor \frac{ \splitDisc  -2}{2}\rfloor}$ on the sub-manifold and combine reduced parameters as $\allParams := [\allParams_\potParams,\allParams_\kinParams] \in \mathbb R^{\splitDisc -2}$. For $\splitDisc$ even, both sub-parameterisations $\allParams_\potParams, \allParams_\kinParams$ are of equal size, i.e.\ $|\allParams_\potParams| = |\allParams_\kinParams| = \frac{\splitDisc - 2}{2}$. However, if $\splitDisc$ is odd, the sub-parametrisations cannot be of equal size. The symbolic zero $\beta_\splitDisc=0$ suggests that $\allParams_\kinParams$ has fewer degrees of freedom and hence $\allParams_\kinParams$ has a smaller sub-parameterisation than $\allParams_\potParams$, i.e.\ 
   $|\allParams_\potParams| = \frac{\splitDisc - 1}{2}$ and $|\allParams_\kinParams| = \frac{\splitDisc - 3}{2}$.
The corresponding linear parameter transform $\paramTrans : \R^{\splitDisc - 2} \to \R^{2\splitDisc}$ which will be used in the following can be written in matrix form as,
\begin{equation} \label{eq:matParamTrans}
    [\potParams,\kinParams]^\top = \paramTrans([\allParams_\potParams,\allParams_\kinParams])=
    \begin{pmatrix}
        A & 0 \\
        0 & B 
    \end{pmatrix}
    \begin{bmatrix}
        \allParams_\potParams \\
        \allParams_\kinParams 
    \end{bmatrix}
    +
    \begin{bmatrix}
        C \\
        D 
    \end{bmatrix}
    \, \Leftrightarrow \,
    A \allParams_\potParams + C = \potParams, \,
    B \allParams_\kinParams + D = \kinParams.
\end{equation}
Explicit expressions for the  matrices $A \in \mathbb R^{\splitDisc\times \lfloor \frac{ \splitDisc  -1}{2}\rfloor}$, $B\in \mathbb R^{\splitDisc\times \lfloor \frac{ \splitDisc  -2}{2}\rfloor}$  and the vectors $C, D\in \mathbb R^\splitDisc$ can be found in Appendix~\ref{sect:paramTrans}. The loss function in \eqref{eq:lossFn} and optimisation problem in \eqref{eq:lossArgMin} becomes a function of the parameters $\allParams \in \mathbb R^{\splitDisc -2}$, namely
\begin{equation} \label{eq:transLossFn}
    \mlLossFn(\allParams) = \E{\numState{0} \sim \stateVarDist} \left\lVert \numFlowComp{\difEq}{\totTime}{\stepSize}(\numState{0} ; \paramTrans(\allParams)) - \numState{\text{Ref}}
    \right\rVert_2^2 
    \quad \text{and} \quad 
    \allParams^* = \A{\allParams} \mlLossFn(\allParams).
\end{equation}
In summary, restricting learning to a sub-manifold (1) reduces the dimension of the search space, (2) improves conditioning and (3) guarantees that the learned splitting method is provably consistent.
\subsection{Training algorithm}\label{sec:training_algorithm}
Having set out the machine learning problem, we now discuss our training algorithm for minimising the loss function \eqref{eq:transLossFn}. It is natural to expect that the coefficients of known classical splittings, such as Trotter in \eqref{eq:trottSolve} and Strang in \eqref{eq:strangSolve}, are local minima of the loss function. For higher-dimensional coefficient spaces, local minima might also correspond to compositions of these methods, as discussed in Section~\ref{sec:composing_splitting_methods}. Our goal is to beat preexisting numerical methods for our choice of initial conditions and step sizes, thus care must be taken to avoid poor local minima that may correspond to pre-existing methods.

Although they are usually the standard choice in machine learning applications, pure gradient based methods such as Adam \cite{kingma2014a} risk getting stuck in such local minima which are more prevalent in low- to moderate-dimensional parameter spaces that we consider here. For this reason, we opt to incorporate a global optimisation aspect which attempts to cover a significant fraction of the parameter space with the aim to identify the minima with very low loss values. As described in Algorithm \ref{alg:pipeline}, the key idea is to create a set of candidate parameter values that are distributed widely over the search space, screen these candidates with a simple heuristic to reduce them to a manageable number and then run a SO algorithm with the candidates as starting points to fine-tune the parameter values. We used Adam in all numerical experiments (see Appendix~\ref{sec:optimiser_comparison} for a discussion of other first- and second-order stochastic optimisers).

Let $\mathcal{S}_\text{train} = \{(\numState{0}^{(j)}, \numState{\text{Ref}}^{(j)})\}_{j=1}^{M_{\text{train}}}$ with $u_0^{(j)}\sim \stateVarDist$ be the training dataset with $M_\text{train}$ data points. A validation dataset $\mathcal{S}_\text{valid}$ with $M_\text{valid}$ samples is constructed in the same way. For each subset $\mathcal{B} \subseteq \mathcal{S}_{\text{train}}$ (or equivalently $\mathcal{B} \subseteq \mathcal{S}_{\text{valid}}$), we define the loss function
\begin{equation} \label{eq:transLossFnFinite}
    \mlLossFn(\allParams; \mathcal{B}) = 
    \frac{1}{ \left|\mathcal{B}\right|} \sum_{(\numState{0}, \numState{\text{Ref}}) \in \mathcal{B}}
    \left\lVert \numFlowComp{\difEq}{\totTime}{\stepSize}(\numState{0}; \paramTrans(\allParams)) - \numState{\text{Ref}} \right\rVert_2^2
\end{equation}
which approximates the ``true'' loss in~\eqref{eq:transLossFn} with a finite sample of size $|\mathcal{B}|$.

\begin{algorithm} 
    \caption{Learning Pipeline for (approximately) minimising $\mlLossFn(\allParams)$}
    \label{alg:pipeline}
    \begin{algorithmic}[1]
    \STATE {Generate candidate coefficients $\Gamma =\{\allParams_1,\allParams_2,\dots\}$, by either choosing $\allParams_i\in\mathbb{R}^{\splitDisc -2}$ randomly or as the vertices of a regular grid that covers the bounded domain $\Omega\subset\mathbb{R}^{\splitDisc -2}$ in which we expect to find the global minimum.}
    \FOR{all candidates $\allParams_i\in\Gamma$}
    \STATE {Compute loss $\ell_i := \mlLossFn(\allParams_i; \mathcal{S}_{\text{valid}})$.}
    \ENDFOR
    \STATE {Remove all candidates $\allParams_{j}$ from $\Gamma$ for which $\ell_j>\min_i \{\ell_i\}+\epsilon$ for some $\epsilon>0$.}
    \STATE {Remove all candidates $\allParams_{j}$ from $\Gamma$ for which $\exists \; i$ s.t. $\| \allParams_{j} - \allParams_{i} \| \leq \delta$ and $\ell_j > \ell_i$ for some minimum distance $\delta>0$.}
    \FOR{all candidates $\allParams_i\in\Gamma$}
    \STATE {Improve $\allParams_i$ by applying a fixed number of batched SO steps; for this use different randomly chosen batches $\mathcal{B} \subset \mathcal{S}_{\text{train}}$ at each SO step.}
    \ENDFOR\RETURN  $\allParams_{\min}=\argmin_{\allParams_i\in\Gamma}\{\mlLossFn(\allParams_i; \mathcal{S}_{\text{valid}})\}$.
    \end{algorithmic}
\end{algorithm}
Observe that the variance of the estimator in ~\eqref{eq:transLossFnFinite} decreases $\propto  |\mathcal{B}|^{-1}$. As a consequence, for small $|\mathcal{B}|$ the value of $\mlLossFn(\allParams;\mathcal{B})$ for fixed $\allParams$ can vary significantly between different batches $\mathcal{B}$. This is not surprising since the learned splitting method will perform differently for different initial conditions. Hence, when making comparative judgments between different parameters, such as in lines 3 and 10 of Algorithm~\ref{alg:pipeline}, the loss function $\mlLossFn(\allParams_i; \, \cdot \,)$ should be evaluated with a single fixed batch, chosen to be $\mathcal{S}_{\text{valid}}$, for all parameters $\allParams_i$.

Each minimum of the loss function will have some basin of attraction under SO. Because of this, it would be computationally inefficient to fine-tune all candidates identified in the first step of Algorithm~\ref{alg:pipeline} with expensive SO iterations. We therefore only pursue candidates that have a low loss (compared to all other candidates) and which are separated by some distance $\delta>0$. This distance is chosen heuristically to be smaller than the average diameter of all basins of attraction.
While alternative techniques such as particle swarm methods~\cite{noel2010a, horst1989a} have been suggested in the literature, we find that empirically Algorithm~\ref{alg:pipeline} produces good results.

Because second-order convergence, in the limit $\stepSize\rightarrow 0$, of our learned methods is guaranteed by their construction in Section~\ref{sec:parameter_transform}, the learned method will generalise to other initial conditions outside the training distribution $\stateVarDist$. This is a distinctive advantage compared to naive ML approaches, which cannot be expected to generalise in this sense. 
\section{Numerical results} \label{sec:implementation_and_results}
In this section we numerically show the efficiency of learned splitting methods for a representative model problem.
\subsection{Schrödinger's equation} \label{sec:schrodinger}
To demonstrate the advantages of our approach, we consider the one-dimensional Schrödinger equation $\imag \dot{\stateVar}(\spaVar, \timeVar) = [\potOpr(\spaVar)-\lapOpr] \stateVar(\spaVar, \timeVar)$ where $\lapOpr=\partial^2/\partial\spaVar^2$ is the Laplace operator and $\potOpr(x)$ is a real-valued potential which we assume to be monotonically increasing for $|\spaVar|>L$. While typically the problem is defined on $\mathbb R\times [0,\totTime]$, discretisation of the unbounded spatial domain with an equidistant grid would require an infinite number of unknowns and make the solution untractable on a computer with finite memory. Because of this, we restrict the domain to $\Omega\times [0,\totTime]$, where $\Omega=[-L,+L]$ for some $L>0$, and impose periodic boundary conditions in space; these boundary conditions render the Laplace operator diagonal in Fourier space \cite{expsolve}. The solution on $\Omega \times [0,\totTime]$ differs from the one on $\mathbb R \times [0,\totTime]$ by terms that are exponentially suppressed. These differences are small provided $L$ is sufficiently large (to see this, note that the eigenfunction corresponding to energy $E$ can be bounded by $C\exp[-\sqrt{\potOpr(L)-E}]$ for $|\spaVar|>L$ and some constant $C$).

The Schrödinger equation can be seen as an example of a wider class of one-dimensional autonomous PDEs of the form
\begin{equation} \label{eq:oneDimDiffEq}
    \dot{\stateVar}(\spaVar, \timeVar) = \mathcal{F}(\stateVar(\spaVar, \timeVar)), \quad \trueState{\spaVar, 0} = \numState{0}(\spaVar), \quad \trueState{\spaVar, \timeVar} \in \C, \quad \spaVar \in \Omega \subset \R, \quad \timeVar \in [0,\totTime]
\end{equation}
with suitable boundary conditions on $\partial\Omega$. Here, $\mathcal{F}$ (which is $i[\lapOpr - \potOpr(\spaVar)]$ for the Schrödinger equation), is the differential operator acting on the solution $\stateVar(\spaVar, \timeVar)$. The method of lines is used to convert the problem \eqref{eq:oneDimDiffEq} into a system of coupled ODEs of the form \eqref{eq:diffEq} which can be solved numerically. For this, we pick a discretisation of the spatial domain $\Omega$ defined by a set of $\spaDisc$ equally spaced points $\{\spaVar_\spaInd\}_{\spaInd=1}^{\spaDisc}$, \mbox{$\spaVar_\spaInd=(\frac{2\spaInd-1}{\spaDisc}-1)L\in\Omega$}, discretise the differential operator $\mathcal{F}$ and solve for the time-dependent solution vector $\stateVar(\timeVar)=(u_1(t),\ldots,u_M(t))\in\C^\spaDisc$ with $\stateVar_\spaInd(\timeVar) \approx \stateVar(\spaVar_\spaInd, \timeVar)$ for $m=1,\ldots, \spaDisc$. 
For the Schrödinger equation, this allows us to write the problem in the form of the IVP in  \eqref{eq:diffEq} as
\begin{equation} \label{eq:potKinSchroedingerPDE}
    \dot{\stateVar}(\timeVar)  = i[\lap - \pot] \stateVar(\timeVar),\quad\trueState{0} = \numState{0}, \quad \trueState{\timeVar} \in \C^\spaDisc, \quad \timeVar \in [0,\totTime]
\end{equation}
for some initial condition $\numState{0}\in \mathbb C^M$. The Laplace operator $\lapOpr$, and indeed the entire Hamiltonian $\hamil = \potOpr(\spaVar)-\lapOpr$, is self-adjoint. The Hermitian matrix $\lap \in \mathbb C^{\spaDisc\times\spaDisc}$ approximates  $\lapOpr$, and is diagonal in discrete Fourier space. Additionally, $\pot \in \mathbb R^{\spaDisc\times\spaDisc}$ is a diagonal matrix, where the diagonal entries are given by $\pot_{m,m} = \potOpr(\spaVar_m)$. The natural splitting for \eqref{eq:potKinSchroedingerPDE}  is to use $\difEqComp{1} = - \imag \pot$ and $\difEqComp{2} = \imag \lap$ as this lets us exploit the diagonalisability of both these flows. We let $\lap = U^\dagger \lap_{\mathrm{diag}} U$ where $\lapdiag \in \mathbb C^{\spaDisc \times \spaDisc}$ denotes a diagonal matrix and $U \in \mathbb C^{\spaDisc\times\spaDisc}$ is the unitary Fourier transform matrix. Note that applying $U$ and $U^\dagger$ with the Fast Fourier Transformation \cite{cooley1965algorithm} costs $\mathcal{O}(\spaDisc\log(\spaDisc))$ operations. As a consequence, the sub-flows corresponding to $\difEqComp{1}$ and  $\difEqComp{2}$ are given by
\begin{equation} \label{eq:potKinflows}
    \exactFlowComp{1}{\timeVar} = \mathrm{e}^{- \imag \timeVar \pot} \text{ and } \exactFlowComp{2}{\timeVar} = \mathrm{e}^{\imag \timeVar \lap} = U^\dagger \mathrm{e}^{\imag \timeVar \lapdiag} U,
\end{equation}
respectively. The sub-flows can be evaluated  efficiently since $\mathrm{e}^{- \imag \timeVar \pot}$ and $\mathrm{e}^{\imag \timeVar \lapdiag}$ are exponentials of diagonal matrices which can be computed in $ \mathcal{O}(\spaDisc)$ time. Hence, the total cost for one evaluation of $\exactFlowComp{1}{\timeVar}$ and $\exactFlowComp{2}{\timeVar}$ is $\mathcal{O}(\spaDisc)$ and $\mathcal{O}(\spaDisc \log 
\spaDisc)$, respectively. In our numerical experiments we fix $\spaDisc=200$.
\subsection{Implementation}\label{sec:implementation}
Our Python code is publicly available on \href{http://doi.org/10.5281/zenodo.13871934}{Zenodo} and datasets and scripts to recreate all figures are available on \href{https://github.com/hl785/paperFigsRepo}{Github}. All training and inference algorithms were implemented in the JAX library \cite{jax2018github}, which we found to be significantly faster than PyTorch since it allows the very efficient automatic differentiation with respect to the parameters $\potParams,\kinParams$ in products of exponentials that arise from \eqref{eq:splitParam}. Double precision arithmetic was used for all calculations. We  used the Optax package \cite{jaxopt_implicit_diff} for the implementation of Stochastic Optimisation in  line 8 of Algorithm~\ref{alg:pipeline}. The Fast Fourier Transform that is required in the evaluation of the flow $\exactFlowComp{2}{\timeVar} = \mathrm{e}^{\imag \timeVar \lap}$ in \eqref{eq:potKinflows} employs the efficient \texttt{jax.numpy.fft} method
to mimic the implementation in the Expsolve package \cite{expsolve}. During training and validation the analytical flow $\exactFlowComp{\difEq}{\totTime} = \e^{\imag \totTime [\lap - \pot]}$ needs to be computed. For this, the matrix exponential is evaluated by exponentiating the (dense) matrix $\lap - \pot$ with the \verb!scipy.linalg.expm! function, which employs the scaling and squaring algorithm given in \cite{alMohy2010a}. Since this is only required during training and will not have any impact on the performance of the learned timestepping methods, no attempt was made to optimise this operation. 
\subsection{Learning splitting coefficients for Schrödinger's equation}
During training we learn optimal splitting coefficients for the Schrödinger equation given in \eqref{eq:potKinSchroedingerPDE} where $V(x) = x^4 - 10 x^2$ is a quartic function that describes a double-well potential. The minima of the potential are located at $x_\pm = \pm \sqrt{5}$ in this case. The parameter transform described in Section \ref{sec:parameter_transform} is used to enforce consistency and symmetry of the learned splitting methods. 
\subsubsection{Data generation} \label{sec:dataGen}
The training dataset $\mathcal{S}_\text{train} = \{(\numState{0}^{(j)}, \numState{\text{Ref}}^{(j)})\}_{j=1}^{M_{\text{train}}}$ consists of initial conditions and reference solutions, where $\numState{0}^{(j)}$ and $\numState{\text{Ref}}^{(j)}$ are functions evaluated on the spatial grid $\{\spaVar_\spaInd\}_{\spaInd=1}^{\spaDisc}$ introduced in Section \ref{sec:schrodinger}. The reference solutions are calculated using the matrix exponential of the initial conditions. While we have not specified the distribution $\stateVarDist$ of initial conditions explicitly, this distribution is defined implicitly by a sampling algorithm that is controlled by three parameters:  the centre mean $x_{\text{cent}}$, the centre standard deviation $x_{\text{stdDev}}$, and the basis standard deviation $\sigma$. We let $x_{\text{cent}} = -\sqrt{5}$, $x_{\text{stdDev}}=0.1$, and $\sigma=0.5$ for training.

We recursively generate $u_0^{(j)}$ by starting from the previous reference solution $\numState{\text{Ref}}^{(j-1)}$. With a moderate probability we perturb $\numState{\text{Ref}}^{(j-1)}$ by adding the discrete Gaussian $\mathcal{N}(x_0^{(j)}, \sigma)$ then renormalising and/or a random phase shift, or with a low probability we set $\numState{0}^{(j)}$ to be the discrete Gaussian $\mathcal{N}(x_0^{(j)}, \sigma)$, where  $x_0^{(j)} \sim \mathcal{N}(x_{\text{cent}}, x_{\text{stdDev}})$. As $x_{\text{cent}} = x_-$, these discrete Gaussians are clustered around the centre of the left well of the potential $V$. For $j=0$, we proceed as above but start with the discrete Gaussian $\mathcal{N}(\bar x_0,\sigma)$ for some $\bar x_0\sim \mathcal N (x_{\text{cent}}, x_{\text{stdDev}})$ instead of the previous reference solution. The recursive augmentation method is written down explicitly in Algorithm~\ref{alg:training_data_generation}, along with some representative initial conditions $\numState{0}^{(j)}$ in Figure~\ref{fig:sampleInitConds}, in Appendix~\ref{sect:training_data}.

A fixed validation dataset $\mathcal{S}_\text{valid}$ is constructed in the same way. During training, we fix $\totTime = 10$ and $\stepSize = \frac{1}{7}$ (which implies $\timeDisc = 70$), and use the approximate loss function in \eqref{eq:transLossFnFinite}. We find that training with a fixed step size results in learned splitting methods which generalise to other timestep sizes $\stepSize$ and final times $\totTime$ (see Section~\ref{sec:generalised_setups}).
\subsubsection{Training procedure}\label{sec:training_procedure}
\begin{figure}[tb]
    \centering
    \includegraphics[width=0.85\textwidth]{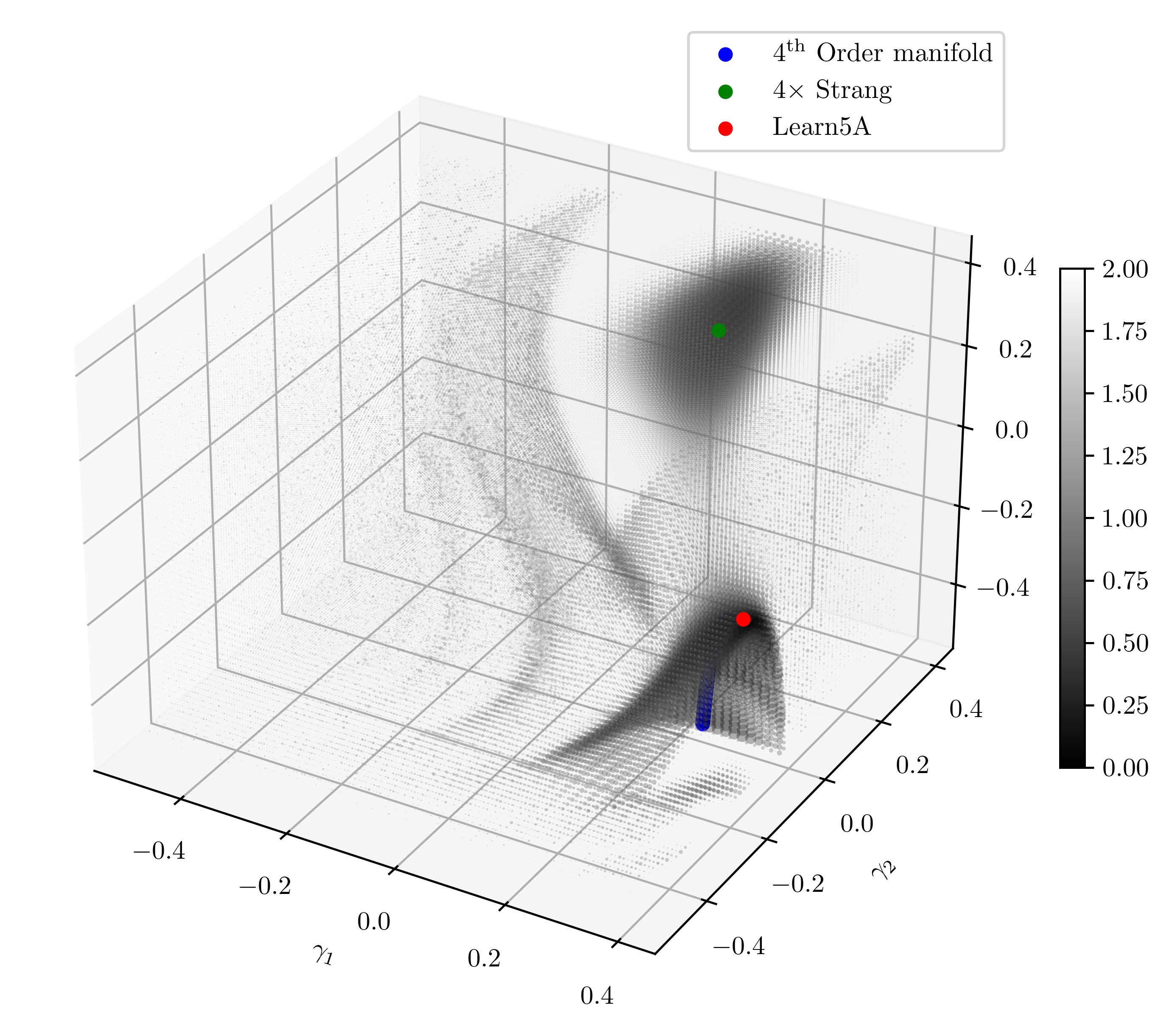}
    \caption{Plot of the loss function $\mlLossFn(\allParams)$ in \eqref{eq:transLossFnFinite} for the Schrödinger equation in \eqref{eq:potKinSchroedingerPDE} with $\totTime=10$ and $\stepSize=\frac{1}{7}$. A section of the one-dimensional manifold of fourth-order accurate methods can be seen in the lower right corner. To obtain this figure, the loss function was evaluated on the fixed validation set $\mathcal{B}=\mathcal{S}_\text{valid}$ with 200 members and for a uniform grid for the values of $\allParams$ in the box $[-0.5,0.4]^3$. Lower loss values are plotted as darker and larger points. 
    \label{fig:lossLandscape}}
\end{figure}

The landscape defined by the loss function \eqref{eq:transLossFnFinite} for splittings of length $\splitDisc = 5$ and evaluated on a fixed validation set $\mathcal{B}=\mathcal{S}_\text{valid}$ of size 200, has three free parameters $\allParams$ and is visualised in Figure~\ref{fig:lossLandscape}.
As this figure shows, the landscape is clearly non-convex and has multiple minima, as expected. Observe that locations of lower values of the loss function appear to be close to lower-dimensional sub-manifolds that define higher order methods, which is consistent with the discussion in Section~\ref{sec:parameter_transform}. The minimum at $\allParams_{\mathrm{Strang}}=[0.125, 0.25, 0.25]$ is readily identified as the composition of four Strang splittings (resulting in a five-stage method) with a validation loss value of $\mlLossFn(\allParams_{\mathrm{Strang}}) = 0.2917$. However, by using the full training algorithm in Algorithm~\ref{alg:pipeline}, we also identified a novel learned splitting, referred to as Learn5A, with $\allParams_{\mathrm{learned}}=[0.3627, -0.1003, -0.1353]$ and a substantially lower validation loss value of $\mlLossFn(\allParams_{\mathrm{learned}}) = 0.02106$. 

\begin{figure}[tb]
    \centering
    \includegraphics[width=0.85\textwidth]{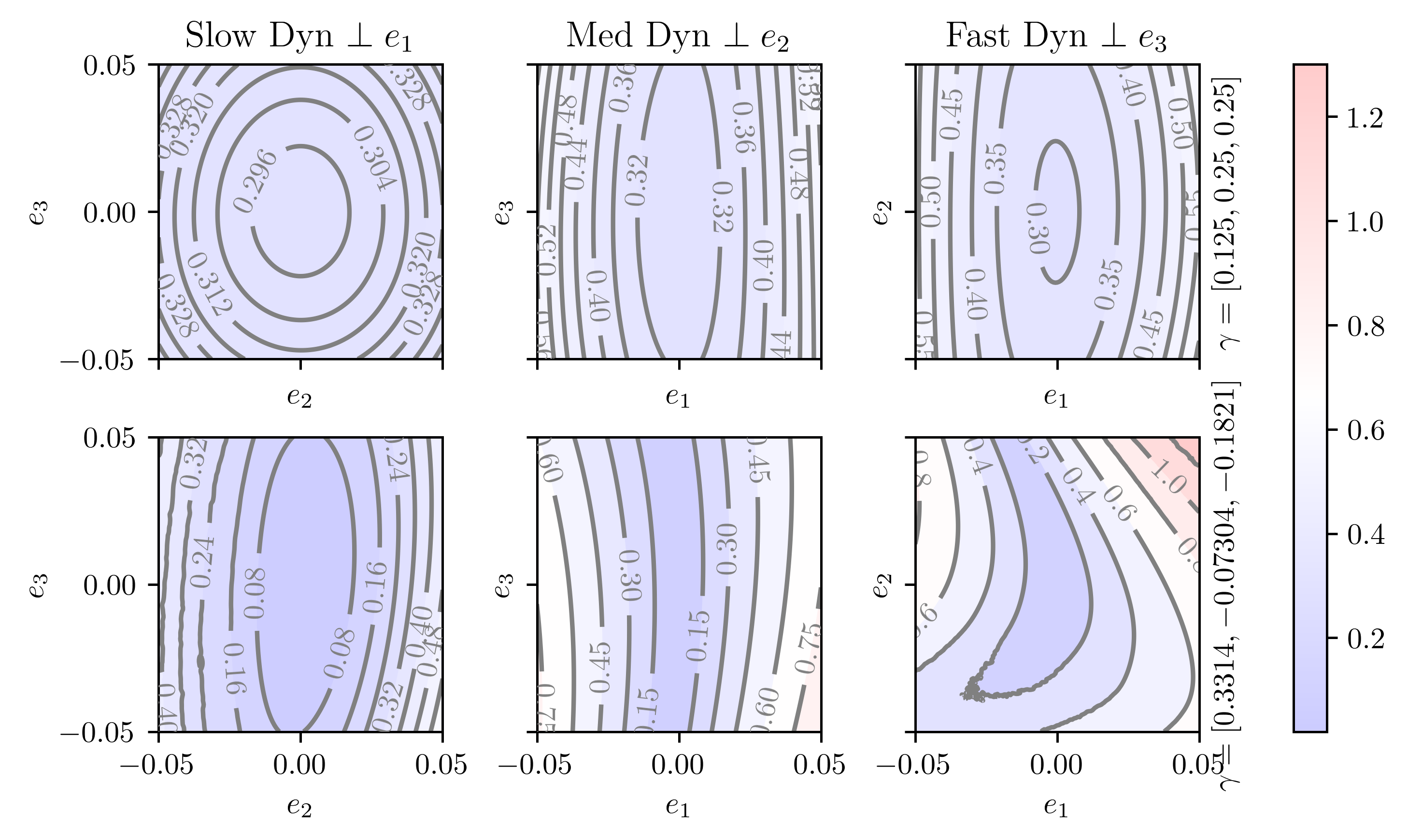}
    \caption{Local environment of the validation loss visualised in Figure~\ref{fig:lossLandscape} around the two minima $\allParams_{\mathrm{Strang}}$ and $\allParams_{\mathrm{valid}}$. The function is plotted in the three planes that are perpendicular to the largest (left), middle (centre), and smallest (right), eigenvalues of the Hessian matrix at the local minima.
    \label{fig:loss2dPlanes}}
\end{figure}

For the training dataset, the landscape will vary between batches, but we expect the qualitative behaviour shown in Figure~\ref{fig:lossLandscape} to be representative for most batches. Due to this stochasticity $\allParams_{\mathrm{learned}}$ will not be a minimum of the loss function when evaluated on $\mathcal{S}_{\text{valid}}$. Instead we consider the nearby validation loss minimum found at $\allParams_{\mathrm{valid}}=[0.3314, -0.07304, -0.1821]$. To gain further insight into the loss function $\mlLossFn(\allParams)$, we plot the validation loss in the proximity of the two minima $\allParams_{\mathrm{Strang}}$ and $\allParams_{\mathrm{valid}}$ in Figure~\ref{fig:loss2dPlanes}. 
The shape of the minima is characterised by their Hessian matrices. The condition numbers of the Hessians at $\allParams_{\mathrm{Strang}}$ and $\allParams_{\mathrm{valid}}$ are 15.5 and 8563.1, respectively. The high condition number of the Hessian matrix at $\allParams_{\mathrm{valid}}$  is consistent with the discussion in Section~\ref{sec:parameter_transform}. This large condition number implies that the naive SGD will be very inefficient, we explored a variety of improved stochastic optimisers (SOs), including both first- and second-order methods, in our training pipeline (line 8 of Algorithm~\ref{alg:pipeline}); further details can be found in Appendix \ref{sec:optimiser_comparison}. However, all results reported in the rest of this section were obtained with the widely used Adam method \cite{kingma2014a} which we find to be efficient and robust. The minimum $\allParams_{\mathrm{learned}}=[0.3627, -0.1003, -0.1353]$ was found with a fixed learning rate of 0.01, where the parameters are initialised with the candidates identified in the initial exploration of the parameter space (lines 1-6 of Algorithm~\ref{alg:pipeline}). The evolution of $\allParams_{\mathrm{learned}}$ during training is shown in Figure~\ref{fig:paramOptim} in Appendix \ref{sec:training_details}.
\subsubsection{Learning longer splitting methods} \label{sec:learning_longer_splittings}
We also employed the training pipeline described in Section~\ref{sec:training_algorithm} to find longer learned splittings with $\splitDisc=8$ stages and six free parameters. For this, we generated 75,000 random parameter candidates, where we implicitly ensured that some of the parameter values are negative and hence potentially  near methods higher than second order. We then evaluated the loss on a consistent validation set (for fairness of comparison), selected the 100 parameter candidates with the lowest losses and removed all candidates that were within an (Euclidean) distance of less than 0.75 to parameter candidates with lower loss values. This resulted in nine candidates that were explored further. For these, the Adam optimisation step in line 8 of Algorithm~\ref{alg:pipeline} was performed for 250 iterations with a learning rate schedule starting at $0.02$ and decreasing exponentially with a decay rate of $0.995$. As before, the parameter transform $\paramTrans(\allParams)$ in \eqref{eq:matParamTrans} was used to ensure that the learned methods are consistent and symmetric. Two of the nine parameter candidates  converged to unknown splittings with low losses, which we refer to as Learn8A and Learn8B. Along with Learn5A found in Section \ref{sec:training_procedure}, the splitting coefficients $\allParams$ after training are given in Table~\ref{tab:learned_coefficients} and their evolution during training is shown in Figure~\ref{fig:paramOptim} in Appendix~\ref{sec:training_details}. 
\begin{table}[ht]
\begin{center}
    \begin{tabular}{ |c|c|c| } 
        \hline
        Splitting & $\splitDisc$ & $\allParams$ \\ 
        \hline \hline
        Trotter & 1 & --- \\ 
        \hline        
        Strang & 2 & [] \\ 
        \hline
        Yoshida & 4 & $[0.67560, 1.35120]$ \\ 
        \hline
        $4 \times$ Strang & 5 & $[0.125, 0.25, 0.25]$ \\ 
        \hline
        Learn5A & 5 & $[0.3627, -0.1003, -0.1353]$ \\ 
        \hline
        Learn8A & 8 & $[
             0.2135,
            -0.0582,
             0.4125,
            -0.1352,
             0.4443,
            -0.0251
        ]$ \\ 
        \hline
        Learn8B & 8 & $[
             0.1178,
             0.3876,
             0.3660,
             0.2922,
             0.0564,
            -0.0212
        ]$ \\ 
        \hline
    \end{tabular}
    \caption{Coefficients $\allParams$ of existing and learned splitting methods with different numbers of stages $\splitDisc$. The non-symmetric Trotter given in \eqref{eq:trottSolve} can not be symmetricly parameterised. The Strang method, given in \eqref{eq:strangSolve}, is completely defined by symmetry and consistency. \label{tab:learned_coefficients}}
\end{center}
\end{table}

\begin{figure}[tb]
    \centering
    \includegraphics[width=0.85\textwidth]{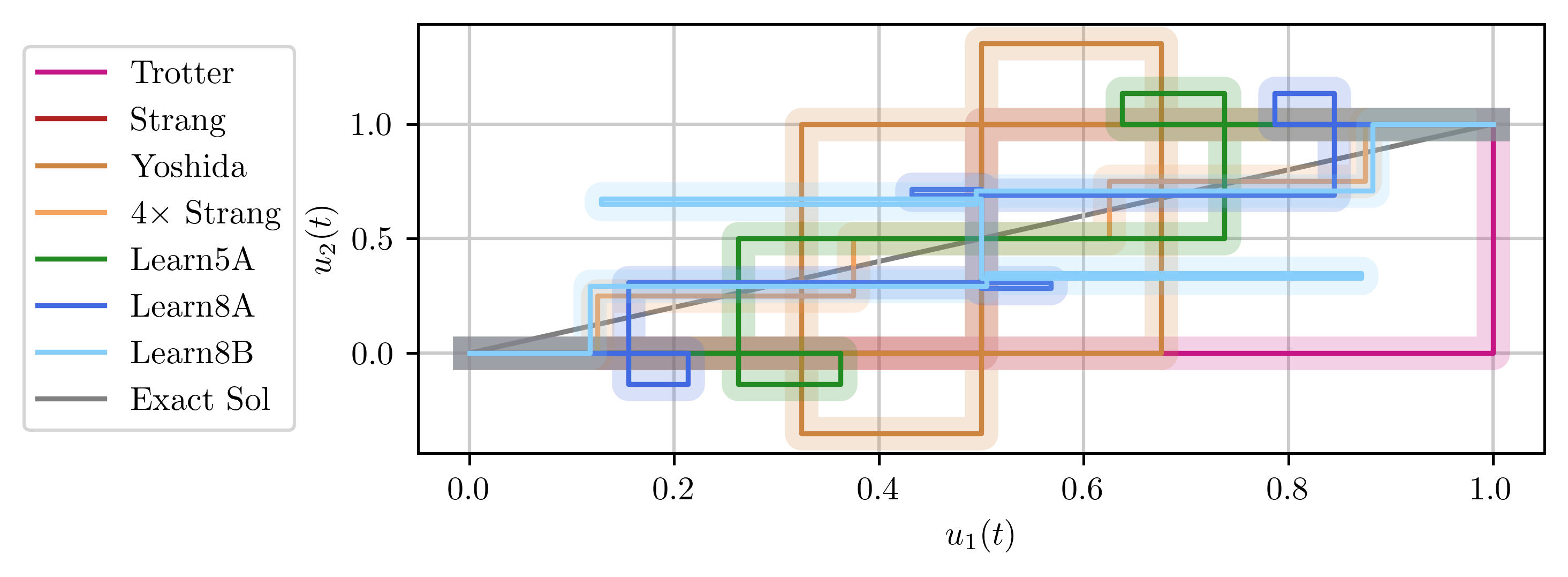}
    \caption{Comparison of the paths that visualise the fluxes \eqref{eqn:tau_flux} for the IVP $\dot{\stateVar}(\timeVar) = [1,1]^\top$, $\trueState{0} = [0,0]^\top$, $\timeVar \in [0,1]$. Results are shown for known numerical methods and our new learned splitting methods.
    \label{fig:splitVisLearn}}
\end{figure}

Naturally, the question arises whether the learned methods are similar to known, existing methods. Since this is difficult to infer from the numerical value of $\allParams$, we  describe a graphical technique for visualising different splitting methods. For this, a method defined by the coefficients $( \potParams_1,\dots,\potParams_\splitDisc,\kinParams_1,\dots,\kinParams_\splitDisc)$ is represented by the continuous curve  obtained by joining line segments of (oriented) lengths $\potParams_1,\kinParams_1,\potParams_2,\kinParams_2,\ldots,$ $\potParams_K,\kinParams_K$ in this order such that the segments associated with $\potParams_j$ are parallel to the horizontal axis while the segments corresponding to $\kinParams_j$ are parallel to the vertical axis. Figure~\ref{fig:splitVisLearn} uses this technique to visualise splitting methods discussed in this work. As explained in Appendix~\ref{sec:splitting_vsualisation}, this construction can be interpreted as the solution of an initial value problem.  As none of the learned paths are covered by another path in their entirety in Figure~\ref{fig:splitVisLearn}, we conclude that the learned methods are novel since they differ from existing methods. We expect that the associated splitting methods also behave fundamentally different for more complicated IVPs, in particular the Schrödinger equation considered here.

\begin{figure}[tb]
    \centering
    \includegraphics[width=0.85\textwidth]{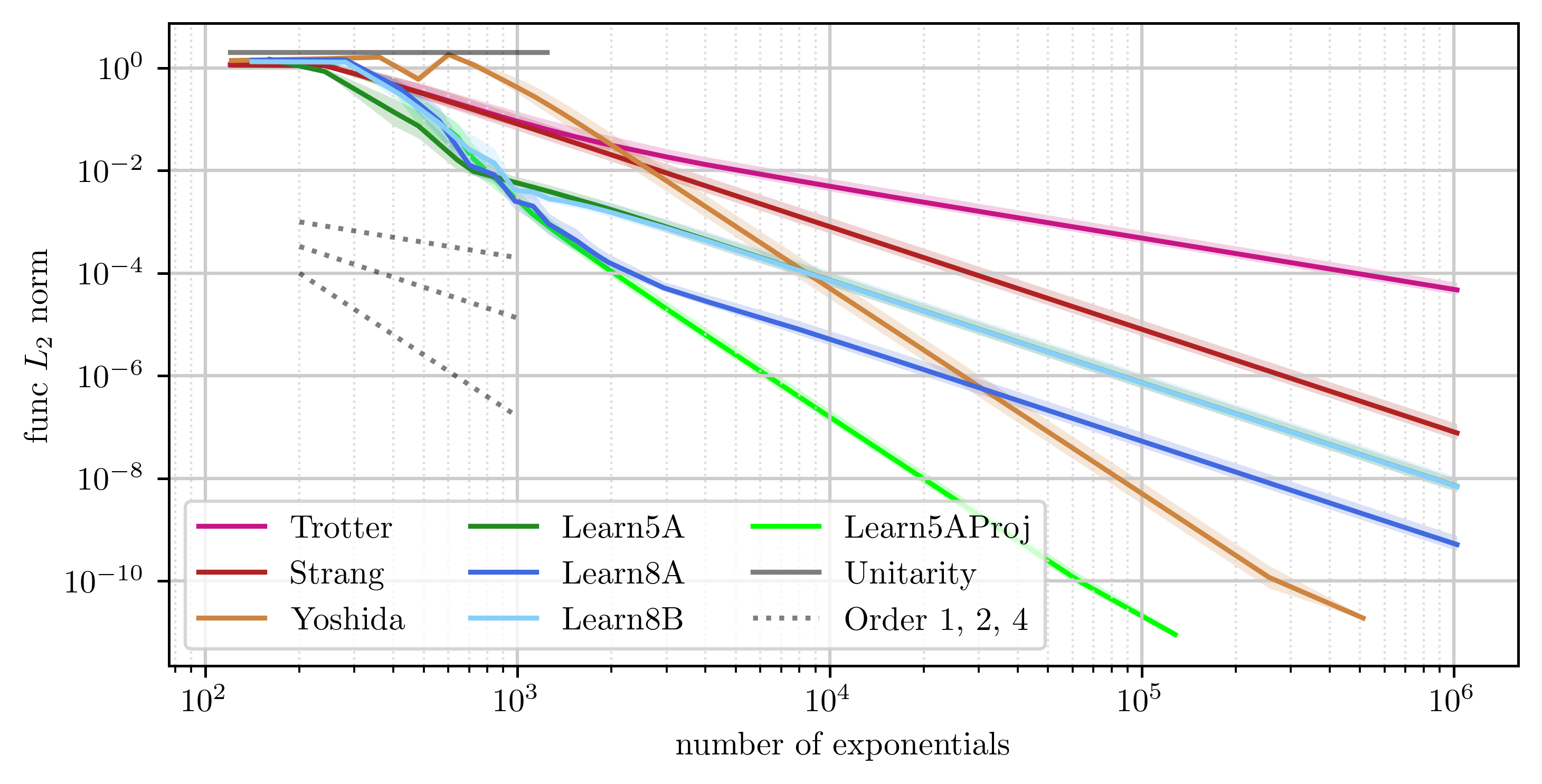}
    \caption{Average $L_2$-error $\left\lVert \numFlowComp{\difEq}{10}{\stepSize}(\numState{0}) -\exactFlowComp{\difEq}{10}(\numState{0}) \right\rVert_2$ at the final time $\totTime=10$ on a validation set of size 200. The median is shown with a line, with a shaded envelope indicating the 15.9\% and 84.1\% quantile. The solid horizontal black vertical line corresponds to the unitarity bound which arises form the fact that the $L_2$-norm of the numerical flow cannot exceed 1. Learn5AProj is introduced in Section \ref{sec:order_of_learned_methods}.
    \label{fig:lossConv}}
\end{figure}

Next, we assess the performance and efficency of our learned methods by comparing them to three well-known splittings, namely Trotter, Strang and Yoshida. As can be seen from Table \ref{tab:learned_coefficients} and Figure~\ref{fig:splitVisLearn}, the different methods involve different numbers of stages $\splitDisc$. This makes comparisons harder, as comparing Strang with $4 \times$ Strang is not fair unless their different cost is taken into account. To address this, we plot the error against the number of exponential evaluations, i.e.\ the number of sub-flow evaluations. Since the number of sub-flow evaluations can also be regarded as a proxy for the computational cost, this allows a fair comparison of the methods with different numbers of stages as we can no longer ``improve'' a method by repeating it many times in a single step of a longer method. Figure~\ref{fig:lossConv} shows the expected $L_2$-error at the final times $\totTime=10$ over the validation set as a function of the number of exponential evaluations (which is proportional to the inverse of the timestep size $\stepSize$) for various methods. 
As expected we observe from Figure~\ref{fig:lossConv} that the learned splittings Learn5A, Learn8A and Learn8B converge quadratically for $\stepSize \rightarrow 0$. Interestingly, they show a faster initial decay for larger $\stepSize$ and this is discussed in more detail in Section~\ref{sec:order_of_learned_methods}. 

\begin{figure}[tb]
    \centering
    \includegraphics[width=0.85\textwidth]{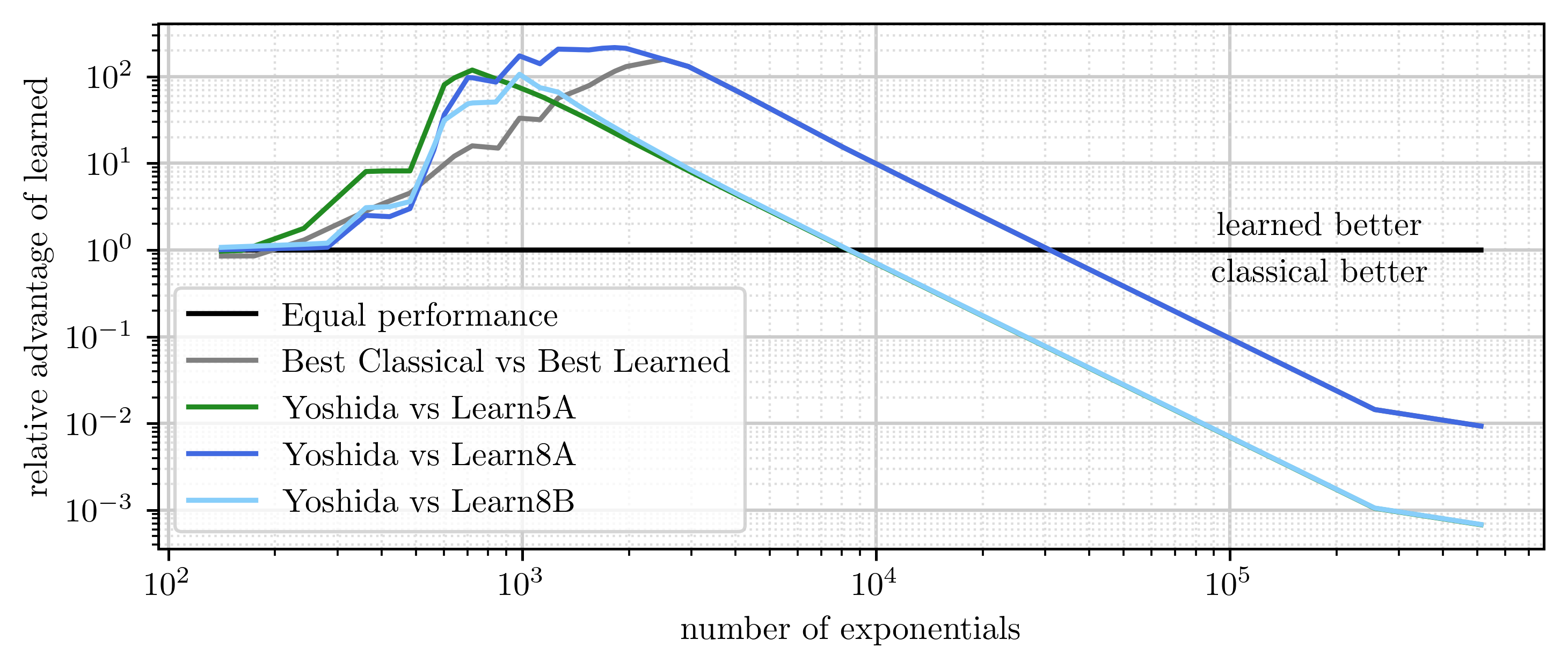}
    \caption{Comparison of the classical vs learned methods relative expected $L_2$-errors, as plotted in Figure~\ref{fig:lossConv}. The plot shows the $L_2$-error of a classical method divided by the $L_2$-error of a learned method. Hence, values larger than one represent the learned method being superior. We define the best classical and learned method as the element-wise minimum of the respective families. 
    \label{fig:lossRelAdv}}
\end{figure}

While the fourth-order Yoshida integrator performs better at high resolutions, Figure~\ref{fig:lossConv} demonstrates that our learned methods outperform all other classical methods for smaller computational budgets: one can either achieve a significant speed-up by using larger step sizes or improve accuracy by orders of magnitude with the same step size. We illustrate this in Figure~\ref{fig:lossRelAdv}, where we plot the relative advantage (in terms of accuracy) in comparison to Yoshida as a function of the number of exponential evaluations. The figure confirms that our learned splitting methods can be up to two orders of magnitude more accurate than Yoshida at the same cost. This is further quantified in Table~\ref{tab:learned_advantage} where we estimated the loss values after 2506 subflow evaluations; our learned methods have the highest relative advantage for this timestep size. The final two columns of the table show the relative accuracy of the methods as well as the gain in speed. The latter is defined as the relative decrease in the number of subflow  evaluations compared to Yoshida when both methods reach the same $L_2$ error. We conclude that we have indeed learned splitting methods that outperform classical numerical methods such as Trotter, Strang, and Yoshida for small computational budgets.
\begin{table}[ht]
\begin{center}
    \begin{tabular}{ |c|c|c|c| }  
        \hline
        \multirow{2}{*}{\textbf{Splitting}} & \textbf{$L_2$-error for} & \textbf{Rel accuracy} & \textbf{Rel speed} \\ 
        & \textbf{$2506$ subflows} & \textbf{vs Yoshida}  & \textbf{vs Yoshida} \\
        \hline \hline
        Trotter & 0.023247 & 0.55 & 
        0.84 \\ 
        \hline        
        Strang & 0.012862 & 1.00 &  1.00 \\ 
        \hline
        Yoshida & 0.012864 & 1.00 &  1.00 \\ 
        \hline
        Learn5A & 0.001121 & 11.47 & 1.84 \\ 
        \hline
        Learn8A & 0.000081 & 158.75 & 3.55 \\ 
        \hline
        Learn8B & 0.001029 & 12.50 & 1.88 \\ 
        \hline
    \end{tabular}
    \caption{\label{tab:learned_advantage}Comparison of loss values after 2506 subflow evaluations, corresponding relative accuracy and gain in speed for various splitting methods.}
\end{center}
\end{table}
\subsubsection{Generalisation to other setups}\label{sec:generalised_setups}
Due to the cost associated with training, learning an optimised splitting method is most justified if it generalises to other setups. Sections \ref{sec:training_procedure} and \ref{sec:learning_longer_splittings} demonstrate that our learned methods are convergent and extrapolate from the training set to the unseen validation set drawn from the same distribution, as expected for any sensible ML approach. However, we note that all methods were trained for a fixed final time $\totTime = 10$ and timestep size $\stepSize = \frac{1}{7}$, which corresponds to 561 sub-flow evaluations for Learn5A and 981 sub-flow evaluations for Learn8A and Learn8B. As shown in Figure~\ref{fig:lossRelAdv}, our learned methods maintain and indeed improve their advantages when generalised to other numbers of subflow evaluations. Next, we investigate the generalisation to other final times, different initial conditions and for perturbations of the double-well potential. More specifically, we consider the following setups:
\begin{itemize}
\item\textbf{Final time:} The final time $\totTime$ was set to either $\totTime_1=10$ or $\totTime_2=30$.
\item\textbf{Initial conditions:} Recall that the (distribution of) initial conditions in Section~\ref{sec:dataGen} were controlled by three parameters $x_{\text{cent}}$, $x_{\text{stdDev}}$ and $\sigma$ and we denote the distribution of initial conditions by $\stateVarDist(x_{\text{cent}}, x_{\text{stdDev}}, \sigma)$. Here we vary the (distribution of) initial conditions by either setting $\stateVarDist_1 = \stateVarDist(-\sqrt{5}, 0.1, 0.5)$, $\stateVarDist_2 = \stateVarDist(\sqrt{5}, 0.2, 0.5)$, or $\stateVarDist_3 = \stateVarDist(-\sqrt{15}, 0.05, \sqrt{0.1})$
\item\textbf{Shape of the potential:}
The shape of the double-well potential $\potOpr(x)$ is set to either $\potOpr_1(x) = x^4 - 10x^2$, $\potOpr_2(x) = x^4 - 10x^2 - 10x$, $\potOpr_3(x) = 3x^4 - 50x^2 + 20x$, or $\potOpr_4(x) = x^4 - 30x^2$; note that this includes non-symmetric potentials.
\end{itemize}
\begin{figure}[tb]
    \centering
    \includegraphics[width=0.85\textwidth]{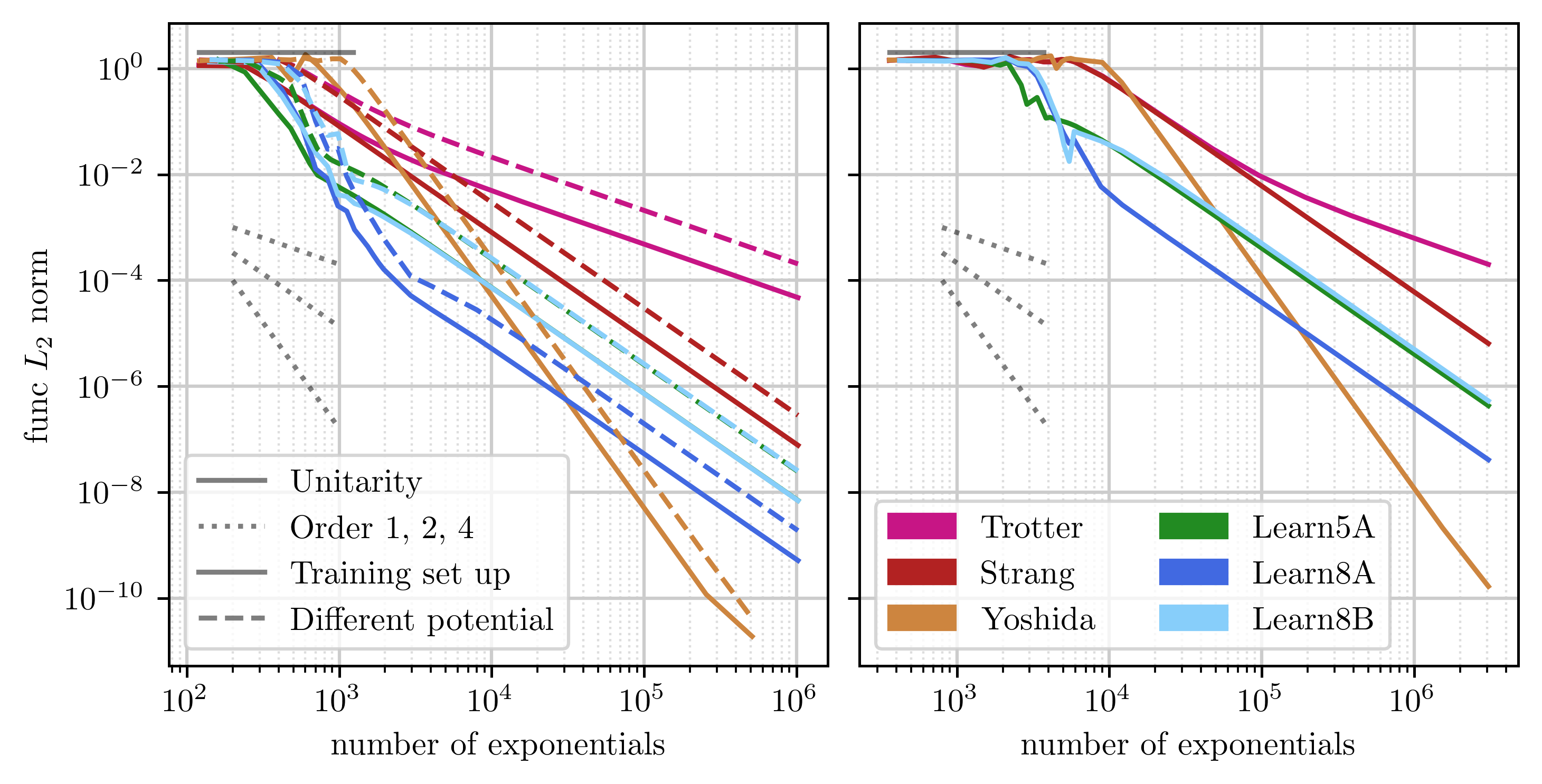}
    \caption{$L_2$ error for generalisations to additional datasets with different final times, initial conditions and potentials. The same quantities as in Figure \ref{fig:lossConv} are plotted. The solid line in the left figure shows the $L_2$-error for the original training data, i.e.\ $\totTime = \totTime_1$, $\stateVarDist = \stateVarDist_1$, and $\potOpr = \potOpr_1$, for reference. The dashed line in the left figure shows $L_2$-error for a changed potential, i.e.\ $\totTime = \totTime_1$, $\stateVarDist = \stateVarDist_1$, and $\potOpr = \potOpr_2$. In the right figure all parameters are changed simultaneously: the solid line  shows the $L_2$-error for different final time $\totTime = \totTime_2$, initial conditions $\stateVarDist = \stateVarDist_2$, and potential $\potOpr = \potOpr_3$. \label{fig:lossConvGen}}
\end{figure}
Figure~\ref{fig:lossConvGen} illustrates the $L_2$-error for different final time, initial conditions and the shape of the potential. Both variations of individual parameters (left) and of all parameters (right) are shown. We omit to plot the $L_2$-error for the case where only the initial conditions are changed, i.e.\ $\totTime = \totTime_1$, $\stateVarDist = \stateVarDist_2$, and $\potOpr = \potOpr_1$, as it is visually indistinguishable from the $L_2$ error on the original training data. We conclude that our learned methods do indeed generalise to IVPs with different final times, initial conditions and potentials and maintain their advantage compared to classical integrators for larger timestep sizes under these variations.

\begin{figure}[tb]
    \centering
    \includegraphics[width=0.47\textwidth]{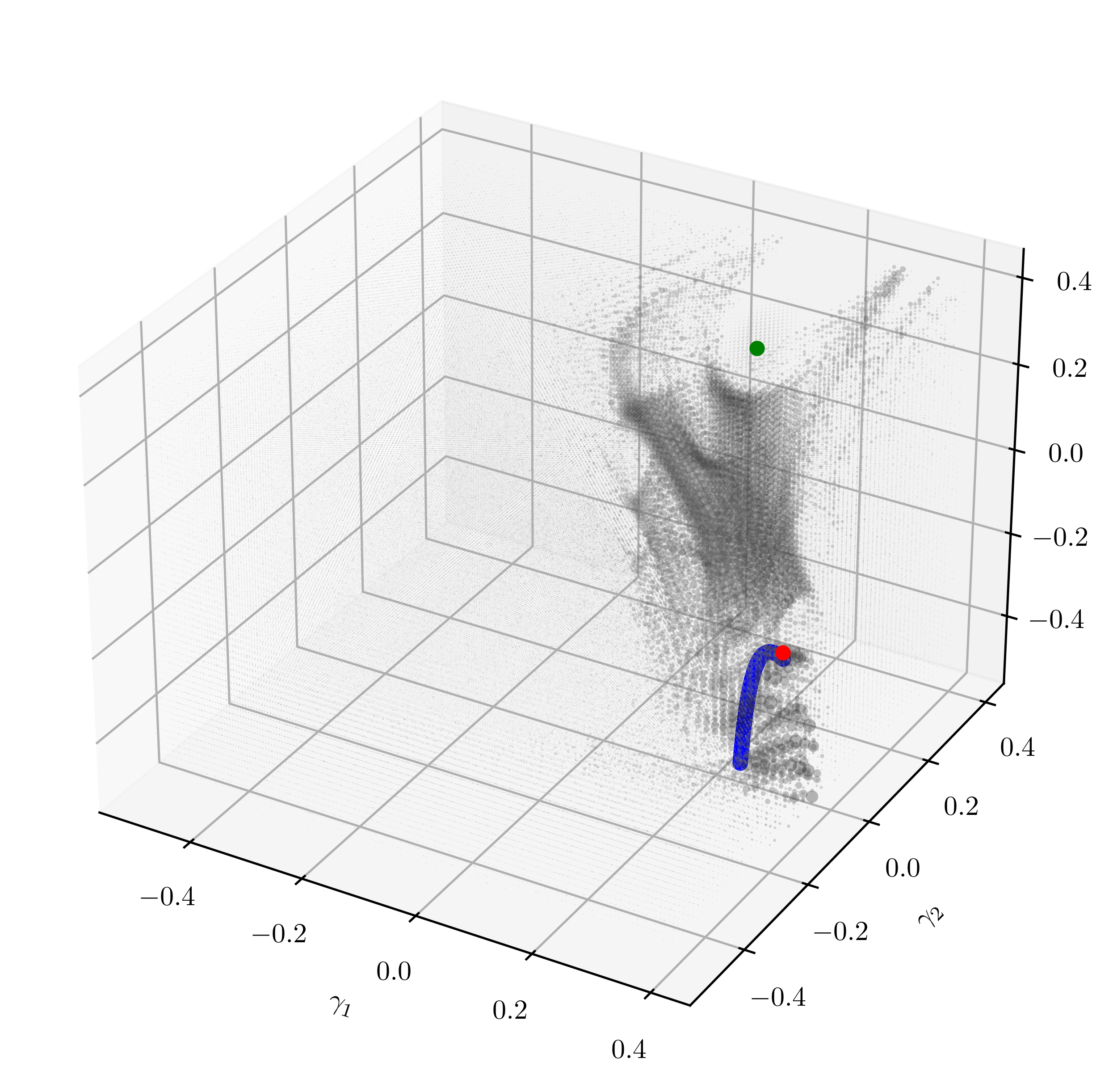}
    \includegraphics[width=0.52\textwidth]{lossLandscape}
    \caption{Left error landscape generated with the distribution of IVPs defined by $\totTime = \totTime_1$, $\stateVarDist = \stateVarDist_3$, and $\potOpr = \potOpr_4$. By comparing with Figure~\ref{fig:lossLandscape}, reproduced right, we clearly see that the landscape, and therefore gradients and minima are problem-dependent. This shows that our pipleline and the learned methods are indeed adapted to the training set of IVP distributions.
    \label{fig:lossLandGen}} 
\end{figure}

To demonstrate that we have indeed adapted to our distribution of IVPs (characterised by different distributions of initial conditions and potentials) and and have not simply learned universal splittings (that we may have been able to find with traditional techniques) we illustrate that the loss landscape is significantly distribution dependant. To see this we consider a distribution of initial conditions with sufficiently deep wells and sufficiently low energy initial conditions so that the training dataset illustrates only negligible quantum tunneling, i.e.\ $\totTime = \totTime_1$, $\stateVarDist = \stateVarDist_3$, and $\potOpr = \potOpr_4$. The associated  landscape is shown in Figure~\ref{fig:lossLandGen} and significantly different from the landscape in Figure~\ref{fig:lossLandscape}. In particular, we note that $4\times$ Strang has ceased to be a local minimum.
\subsubsection{Order of learned methods}\label{sec:order_of_learned_methods}
Recall that for a given method,  $\stepSize$ is inversely proportional to the number of exponentials $\timeDisc$ required to integrate to a fixed final time $\totTime$. With this in mind, Figure~\ref{fig:lossConv} shows that  asymptotically the $L_2$-error of the learned methods is quadratically convergent in the timestep size $\stepSize$. This is of course to be expected since the methods were explicitly constructed to be consistent and symmetric. Interestingly, before reaching this asymptotic convergence rate, there is a significant drop of the error for larger values of $\stepSize$. To quantify this, we consider a Taylor expansion of the $L_2$-error. As the odd terms of the Taylor expansion are symbolically zero due to the parameter transform enforcing symmetry, we expand
\begin{equation}
    E(\stepSize) = C_2 h^2 + C_4 h^4 + C_6h^6 + \mathcal{O}(h^{8}), \label{eqn:error_taylor_expansion}
\end{equation}
for some non-negative constants $C_2,C_4,C_6 \ge 0$. If $0 < C_2 \ll C_4,C_6$, then the methods are formally second order ($C_2\ne 0$), but the error will decrease approximately with a higher power of $\stepSize$ for larger stepsizes where the term $C_2 h^2$ in \eqref{eqn:error_taylor_expansion} is small relative to the higher order terms.
\begin{figure}[tb]
    \centering
    \includegraphics[width=0.85\textwidth]{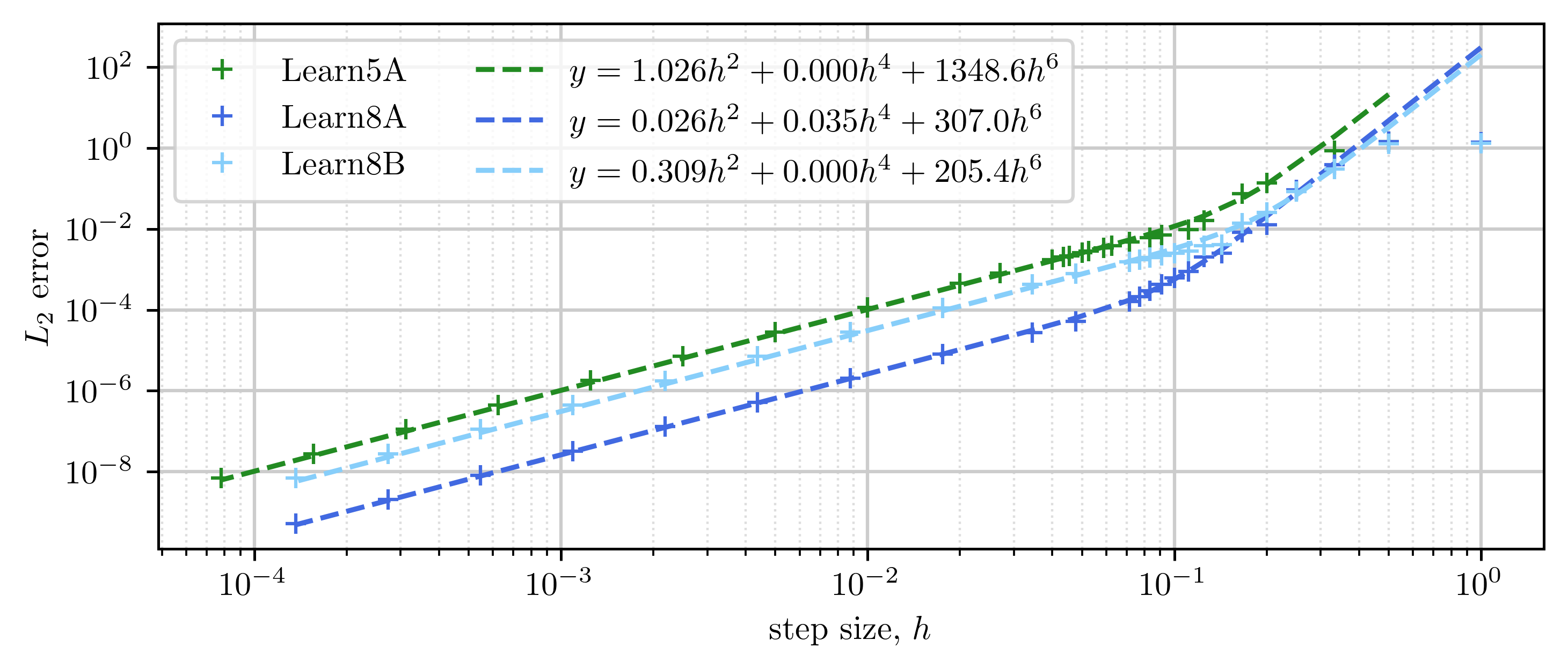}
    \caption{Fit of the coefficients $C_{2j}$ for $j=1,2,3$ in the Taylor expansion \eqref{eqn:error_taylor_expansion} to the data points $(x_i,y_i)$ as a function of the timestep size $\stepSize$. 
    \label{fig:bestFitCoefs}}
\end{figure}
Figure~\ref{fig:bestFitCoefs} shows a fit of the leading terms in the expansion \eqref{eqn:error_taylor_expansion} to the $L_2$-error, excluding the two data points with the largest timesteps where the sixth order polynomial approximation is not appropriate. This fit is based on the assumption that our data is of the form $(x_i, y_i)$ and obeys the model $y_i = \sum_{j = 1}^3 C_{2j} x_i^{2j}$ with coefficients $C_{2j}$. A least-squares fit leads to the minimisation problem $\A{C \in \R^3} \mathcal{L}(C)$ where $C=(C_{2},C_{4},C_{6})\in \mathbb R^3$ and $\mathcal{L}(C) = \sum_{i} (\sum_{j = 1}^3 C_{2j} \phi_j(x_i) - y_i)^2$. To enforce non-negativity of the coefficients, we set $C_{2j} = (\tilde C_{2j})^2$ and fit the logarithm of the error, i.e.\ we find $\A{\tilde C \in \R^3} \mathcal{\tilde L}(\tilde C)$ with $\tilde C=(\tilde{C}_2,\tilde{C}_4,\tilde{C}_6)$ and $\mathcal {\tilde L}(\tilde C)=\sum_{i} (\log(\sum_{j = 1}^3 (\tilde C_{2j})^2 \phi_j(x_i)) - \log(y_i))^2$. The results of this fit in Figure~\ref{fig:bestFitCoefs} show that visually the agreement with the data looks very reasonable for the Learn8A method, while it is less good for the other two methods. This is because the expansion in \eqref{eqn:error_taylor_expansion} does not consist of orthogonal basis functions and hence will be poor for larger values of $\stepSize$ which constrain the terms $C_4$ and $C_6$ as higher order terms will have non-trivial effects. We stress, however, that our goal is not to find the exact values of the fit coefficients but rather to get an indication of the relative sizes of $C_2$, $C_4$ and $C_6$.
\begin{table}[htb]
        \centering
        \begin{tabular}{|c|c|c|c|}
            \hline
            \textbf{method} & $C_{2}$ & $C_{4}$ & $C_{6}$ \\
            \hline\hline
            Learn5A & 1.026 & 0.000 & 1348.6 \\
            Learn8A & 0.026 & 0.035 & 307.0 \\
            Learn8B & 0.309 & 0.000 & 205.4 \\
            \hline
        \end{tabular}
        \caption{Fit coefficients of the $L_2$-error  expansion in \eqref{eqn:error_taylor_expansion}  as a function of the timestep size $\stepSize$ for the three learned methods, based on the fit shown in Figure \ref{fig:bestFitCoefs}.}
    \label{tab:learnedCoefs}
\end{table}
The numerical values of the fit coefficients $C_2,C_4,C_6$  given in Table~\ref{tab:learnedCoefs} show that for all three learned methods we have $C_2 \ll C_6$. This confirms that although our learned methods are formally of order two, they are very close to methods of higher order. 

\begin{figure}[tb]
    \centering
    \includegraphics[width=0.85\textwidth]{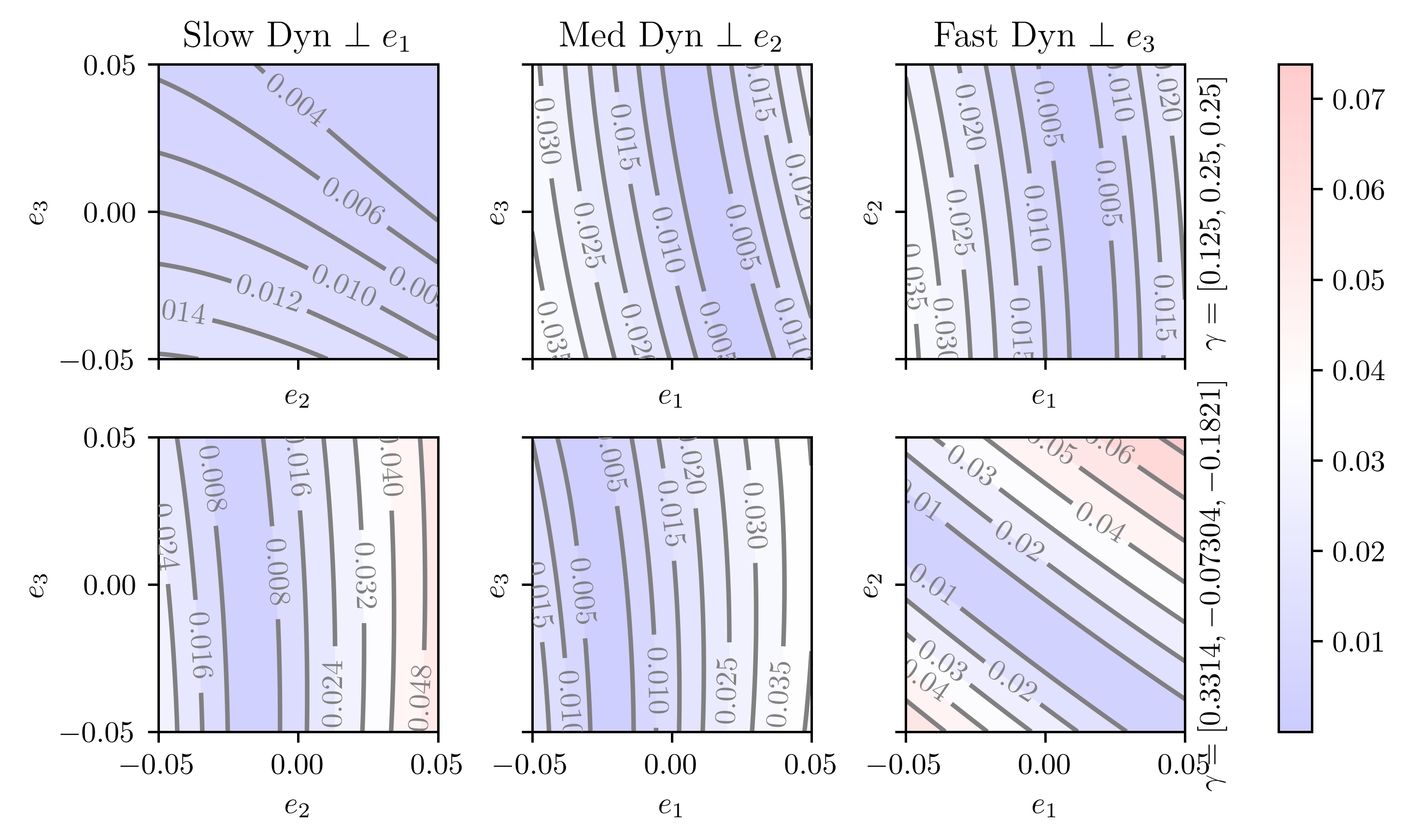}
    \caption{Function $w(\allParams):=|w_{112}|+|w_{122}|$ which acts as a proxy of the violation of the fourth order condition as a function of the parameters $\allParams$. The function is plotted in the vicinity of the minima of the validation loss shown in Figure~\ref{fig:loss2dPlanes}.
    \label{fig:orderConds2dPlanes}}
\end{figure}

Given that empirically the learned methods appear to be ``close to'' fourth order, it is natural to look for higher order methods in the vicinity of our learned methods in the space parametrised by $\allParams$. The key technique for finding such higher order methods is to construct a smooth function of $\allParams$ which vanishes on the manifold that describes higher order methods. As explained in \cite{blanes2024a}, the manifold of methods of at least order four can be implicitly defined with the help of two polynomials $w_{112}, w_{122}$ of the splitting coefficients $(\potParams, \kinParams)$. By using the parameter transform $g(\allParams)$ described in Section \ref{sec:parameter_transform}, we express $w_{112}, w_{122}$ as a function of the independent splitting coefficients $\allParams$. A symmetric method with splitting coefficients $\allParams$ is of order four or higher if $w_{112}(\allParams)=0$ and $w_{122}(\allParams)=0$. For five-stage methods with three free parameters $\allParams$, we can explicitly parametrise the one-dimensional manifold defined by $w_{112}(\gamma)=w_{122}(\gamma)=0$ which contains higher order methods. This manifold is depicted in Figure~\ref{fig:lossLandscape} and Figure~\ref{fig:lossLandGen}. In Figure~\ref{fig:orderConds2dPlanes}, we plot the value of $w(\allParams):=|w_{112}(\allParams)|+|w_{122}(\allParams)|$ as a function of the parameters $\allParams$ in the vicinity of the minima that are  considered in Figure~\ref{fig:loss2dPlanes}. Observe that $w(\allParams)=0$ is an equivalent definition of the higher-order manifold $w_{112}(\allParams)=w_{122}(\allParams)=0$. As  Figure~\ref{fig:orderConds2dPlanes} shows, the manifold of methods of at least order four is close to the minima $\allParams_{\mathrm{valid}}$ of the validation loss, visualised in Figure~\ref{fig:loss2dPlanes}, that was associated with Learn5A. This confirms that our learned methods are indeed close to methods of at least order four. 

Next, we further improve the learned method Learn5A by projecting the coefficient $\allParams_{\mathrm{learned}} = [0.363, -0.100, -0.135]$ to the closest point on the manifold $w(\allParams)=0$. We call this new method, defined by $\allParams_{\mathrm{proj}} = [0.346, -0.112, -0.132]$, Learn5AProj. The $L_2$-error of Learn5A and Learn5AProj is shown in Figure~\ref{fig:lossConv}. As expected the Learn5AProj shows fourth-order convergence. Note, however, that our gradient based method was correct in moving away from the manifold of fourth-order methods as Learn5A does indeed have a smaller error than Learn5AProj for larger timesteps. Overall, this shows that our approach is able to find the most efficient method for a limited, fixed cost budget (i.e. number of exponential evaluations), that these efficient methods may be lower order and lower error constant methods, but also that gradient based methods can approximately recover the classically derived order conditions.
\section{Conclusion and future work}


In this paper, we have developed a machine learning-based framework to find computationally efficient splitting schemes   tailored for limited computational budgets  and adapted to classes of IVPs defined by specific right-hand sides, final times and distributions of initial conditions.  Our approach  maintains  common advantages of classical methods, such as interpretability, generalisability, convergence, and conservation properties, which are usually not guaranteed for machine learning-based approaches. We showed that our framework can learn medium to long splittings with lower errors than classical splittings. It is particularly efficient for large timestep sizes. It should be stressed that our learned methods are novel which was confirmed by comparing them to existing splitting methods. 
There are several ways in which this work can be extended. This includes exploiting the generalisation abilities of a specific learned method. For instance, a learned method could be found by training on a low-dimensional system, which is cheap to simulate, and then applying it to more sophisticated higher-dimensional problems. One could also investigate longer splittings with a larger number of subflows. Further work could also be done to remove the dependence on labelled training data or augment the loss function with regularisers, which penalise the violation of higher-order conditions.
\section*{Acknowledgements}
The computational studies made use of the Nimbus High Performance Computing (HPC) Services at the University of Bath (University of Bath’s Research Computing Group (\url{doi.org/10.15125/b6cd-s854}).

\bibliographystyle{siamplain}
\bibliography{references}


\newpage
\appendix

\section{Triple jump and composition methods} \label{sect:composition}
We describe the triple jump technique \cite[II: 4.2]{hairer2006a} as an example of a composition method. Given a numerical method  $\numFlowCompOneStep{\difEq}{\stepSize}(\, \cdot \, ; \potParams, \kinParams)$ of order $\errorOrder$, we can construct a higher-order numerical method by writing
\begin{equation} \label{eq:tripleJump}
    \numFlowCompOneStep{\difEq}{\stepSize}(\, \cdot \, ; \widetilde{\potParams}, \widetilde{\kinParams})
    =
    \numFlowCompOneStep{\difEq}{\mu_3\stepSize}(\, \cdot \, ; \potParams, \kinParams)
    \circ
    \numFlowCompOneStep{\difEq}{\mu_2\stepSize}(\, \cdot \, ; \potParams, \kinParams)
    \circ
    \numFlowCompOneStep{\difEq}{\mu_1\stepSize}(\, \cdot \, ; \potParams, \kinParams),
\end{equation}
where,
\begin{xalignat}{2}
    \mu_1 &= \mu_3 = \frac{1}{2 - 2^{1/(\errorOrder+1)}},& \mu_2 &= 1 - \mu_1 - \mu_3.
\end{xalignat}
Due to the symmetry of the parameters that define the triple jump, $\mu = [\mu_1, \mu_2, \mu_3]$, within the larger class of composition methods, applying the triple jump to a symmetric method results in a new symmetric method. Hence the new method defined by \eqref{eq:tripleJump} is of order $\errorOrder + 1$ in general and of order $\errorOrder + 2$ if the original method is symmetric. The coefficients $[\widetilde{\potParams}, \widetilde{\kinParams}]$ that define the new method can be constructed from the coefficients $[\potParams, \kinParams]$ of the original method in a straightforward way. 

Applying the triple jump technique to the second order Strang results in a splitting method of order four from the Yoshida family, referred to as Yoshida in this paper. The triple jump technique can be inductively repeated to construct higher-order splitting methods. However, the number of stages grows exponentially in this approach: constructing a method of order $\errorOrder \gg 1$  will result in $\mathcal{O}(3^\errorOrder)$ stages.
\section{Explicit parameter transform} \label{sect:paramTrans}
To define the matrices $A \in \mathbb R^{\splitDisc\times \lfloor \frac{ \splitDisc  -1}{2}\rfloor}$, $B\in \mathbb R^{\splitDisc\times \lfloor \frac{ \splitDisc  -2}{2}\rfloor}$ and the vectors $C, D\in \R^{\splitDisc}$ for the parameter transform in \eqref{eq:matParamTrans}, we distinguish between $\splitDisc$ even and odd. For $\splitDisc$ even, we have
\begin{align*}
    A = 
    \begin{pmatrix}
        I_{\frac{\splitDisc-2}{2}} \\
        -1_{2, \frac{\splitDisc-2}{2} } \\
        J_{\frac{\splitDisc-2}{2} }
    \end{pmatrix}
    \in \R^{\splitDisc \times \frac{\splitDisc-2}{2}}, \qquad 
    & C = 
    \begin{pmatrix}
        0_{\frac{\splitDisc-2}{2}, 1} \\
        0.5_{2, 1} \\
        0_{\frac{\splitDisc-2}{2}, 1}
    \end{pmatrix}
    \in \R^{\splitDisc},
    \\ B = 
    \begin{pmatrix}
        I_{\frac{\splitDisc -2}{2} } \\
        -2_{1, \frac{\splitDisc -2}{2} } \\
        J_{\frac{\splitDisc -2}{2} } \\
        0_{1, \frac{\splitDisc-2}{2} }
    \end{pmatrix}
    \in \R^{\splitDisc \times \frac{\splitDisc-2}{2} }, \qquad
    & D = 
    \begin{pmatrix}
        0_{\frac{\splitDisc-2}{2} , 1} \\
        1 \\
        0_{\frac{\splitDisc-2}{2} , 1} \\
        0
    \end{pmatrix}
    \in \R^{\splitDisc},
\end{align*}
while for $\splitDisc$ odd, we obtain,
\begin{align*}
    A = 
    \begin{pmatrix}
        I_{\frac{\splitDisc - 1}{2}} \\
        -2_{1, \frac{\splitDisc - 1}{2}} \\
        J_{\frac{\splitDisc - 1}{2}}
    \end{pmatrix}
    \in \R^{\splitDisc \times \frac{\splitDisc - 1}{2}}, \qquad
    & C = 
    \begin{pmatrix}
        0_{\frac{\splitDisc - 1}{2}, 1} \\
        1 \\
        0_{\frac{\splitDisc - 1}{2}, 1}
    \end{pmatrix}
    \in \R^{\splitDisc},
    \\ B = 
    \begin{pmatrix}
        I_{\frac{\splitDisc - 3}{2}} \\
        -1_{2, \frac{\splitDisc - 3}{2}} \\
        J_{\frac{\splitDisc - 3}{2}} \\
        0_{1, \frac{\splitDisc - 3}{2}}
    \end{pmatrix}
    \in \R^{\splitDisc \times \frac{\splitDisc - 3}{2}}, \qquad
    & D =
    \begin{pmatrix}
        0_{\frac{\splitDisc - 3}{2}, 1} \\
        0.5_{2, 1} \\
        0_{\frac{\splitDisc - 3}{2}, 1} \\
        0
    \end{pmatrix}
    \in \R^{\splitDisc}.
\end{align*}
Here, $I_s \in \mathbb R^{s\times s}$ denotes the identity matrix and we write $J_s \in \mathbb R^{s\times s}$ for the exchange matrix. Further,  $r_{p, q} \in \mathbb R^{p\times q}$ is filled with the scalar $r$. 
\section{Construction of training datasets} \label{sect:training_data}
\begin{figure}[tb]
    \centering
    \includegraphics[width=0.85\textwidth]{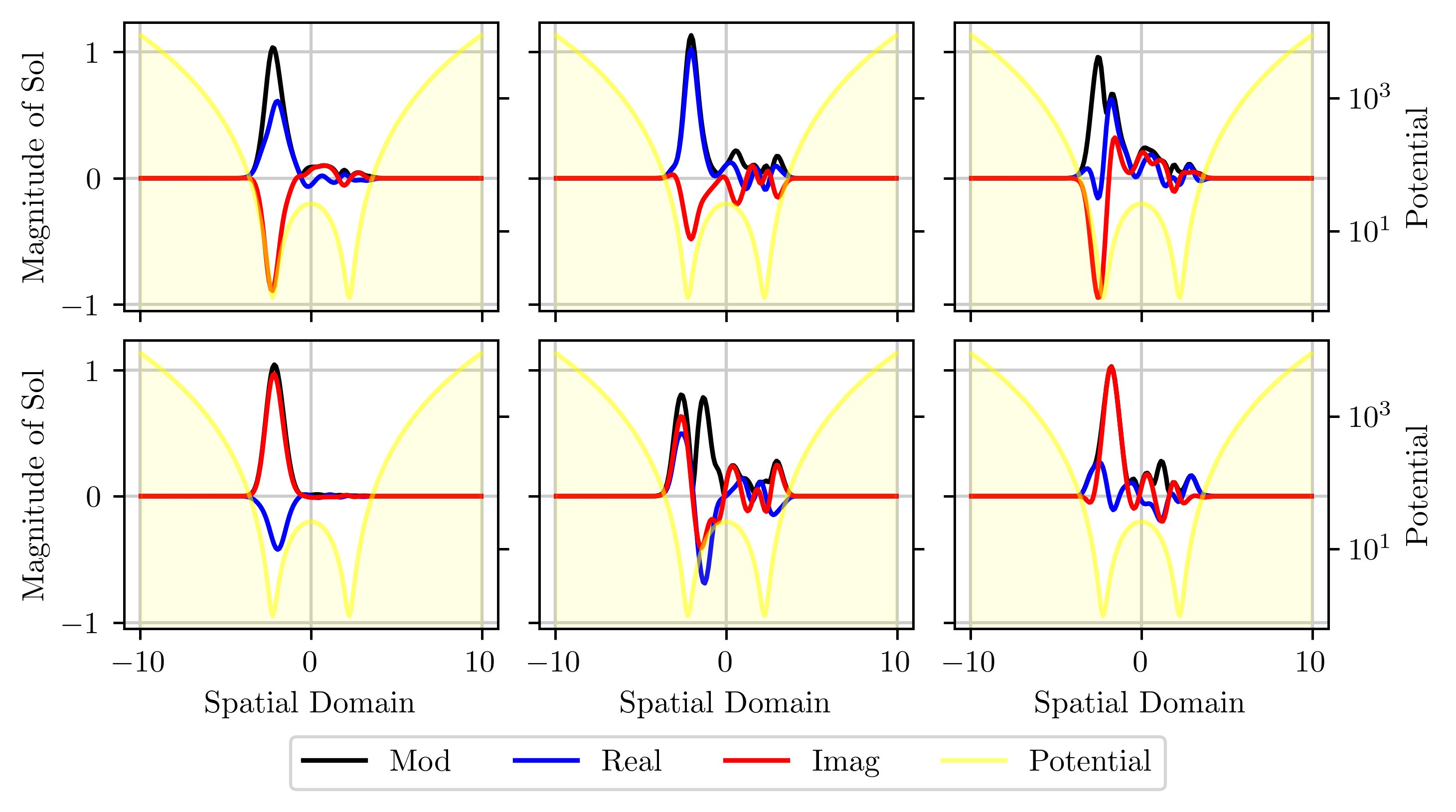}
    \caption{Six randomly chosen initial conditions $u_0^{(j)}$ generated with Algorithm~\ref{alg:training_data_generation}, with modulus, real and imaginary part of $u_0^{(j)}$ depicted in black, blue and red, resptively. In each of the plots, the double-well potential is also shown on a logarithmic scale in the background.
    \label{fig:sampleInitConds}}
\end{figure}

The procedure for generating a single batch $\mathcal{B}$ of random initial conditions which can be used for training and validation is  explicitly described in Algorithm~\ref{alg:training_data_generation}. In lines 5, 8 and 11 of Algorithm~\ref{alg:training_data_generation},  the function $g(\,\cdot\,;\bar x_0,\sigma)$ is a Gaussian with mean $\bar x_0$ and standard deviation $\sigma$, namely
\begin{equation}
    g(x;x_0,\sigma):=\mathcal{Z}\exp\left[-\frac{1}{2}\left(\frac{x-\bar x_0}{\sigma}\right)^2\right].\label{eqn:gaussian_initial_condition}
\end{equation}
Here, $\mathcal{Z}$ is a suitable normalisation constant such that $\sum_{m=1}^{\spaDisc} |g(x_m;x_0,\sigma)|^2=1$, where $x_m,~m=1,\ldots,M$, denotes the spatial discretisation introduced in Section~\ref{sec:schrodinger}. As already remarked in Section \ref{sec:dataGen}, Algorithm \ref{alg:training_data_generation} implicitly defines the distribution $\mathcal{U}$ of initial conditions used for training.
\begin{algorithm} 
    \caption{Generation of training batch $\mathcal{B}=\{u_0^{(0)},u_0^{(1)},\dots,u_0^{(b-1)}\}$.}
    \label{alg:training_data_generation}
    \begin{algorithmic}[1]
        \STATE{Set width $\sigma=0.5$.}
        \FOR {$j=0,1,\ldots,b-1$}
        \STATE{Draw normally distributed random mean $x_0^{(j)} \sim \mathcal{N}(-\sqrt{5},0.1)$}.
        \STATE{Draw uniformly distributed random numbers $\xi_j^{(1)}, \xi_j^{(2)}, \xi_j^{(3)}, \xi_j^{(4)} \sim \mathcal{U}(0,1)$.}
        \STATE{Set $\widetilde{\phi} = 
        \begin{cases}
            g(\,\cdot\,;\overline{x}_0,\sigma) \quad \text{with $\overline{x}_0 \sim \mathcal{N}(-\sqrt{5},0.1)$}, &\text{if $j=0$},\\
            \exactFlowComp{\difEq}{\totTime}(u_0^{(j-1)}), &\text{otherwise}.
        \end{cases}
        $}

        \IF{$\xi_j^{(1)} < 0.5$}
        \STATE{Set $\widetilde{\phi} \mapsto \widetilde{\phi} + g(\, \cdot \,; x_0^{(j)}, \sigma)$ (add Gaussian)}.
        \ENDIF
    
        \IF{$\xi_j^{(2)} < 0.5$}
        \STATE{$\widetilde{\phi} \mapsto \exp[2\pi \imag \xi_j^{(3)}] \widetilde{\phi}$ (apply a random phase shift)}.
        \ENDIF
    
        \IF{$\xi_j^{(4)} < 0.01$}
        \STATE{$\widetilde{\phi} = g(\,\cdot\,;x_0^{(j)},\sigma)$ (reset to Gaussian)}.
        \ENDIF
    
        \STATE{Set $u_0^{(j)} = \widetilde{\phi} / \widetilde{\mathcal{Z}}$ with the normalisation constant $\widetilde{\mathcal{Z}} = \sum_{m=1}^{\spaDisc} |\widetilde{\phi}(x_m)|^2$}.
        \ENDFOR
    \end{algorithmic}
\end{algorithm}
Figure~\ref{fig:sampleInitConds} shows six randomly chosen initial conditions $u_0^{(j)}$ that were generated with Algorithm~\ref{alg:training_data_generation}. Note that the functions $u_0^{(j)}$ are concentrated around the left minimum $x_-=-\sqrt{5}$ of the double-well potential.
\section{Further details on the choice of the stochastic optimisation algorithm}\label{sec:training_details}
\begin{figure}[tb]
    \centering
    \includegraphics[width=0.85\textwidth]{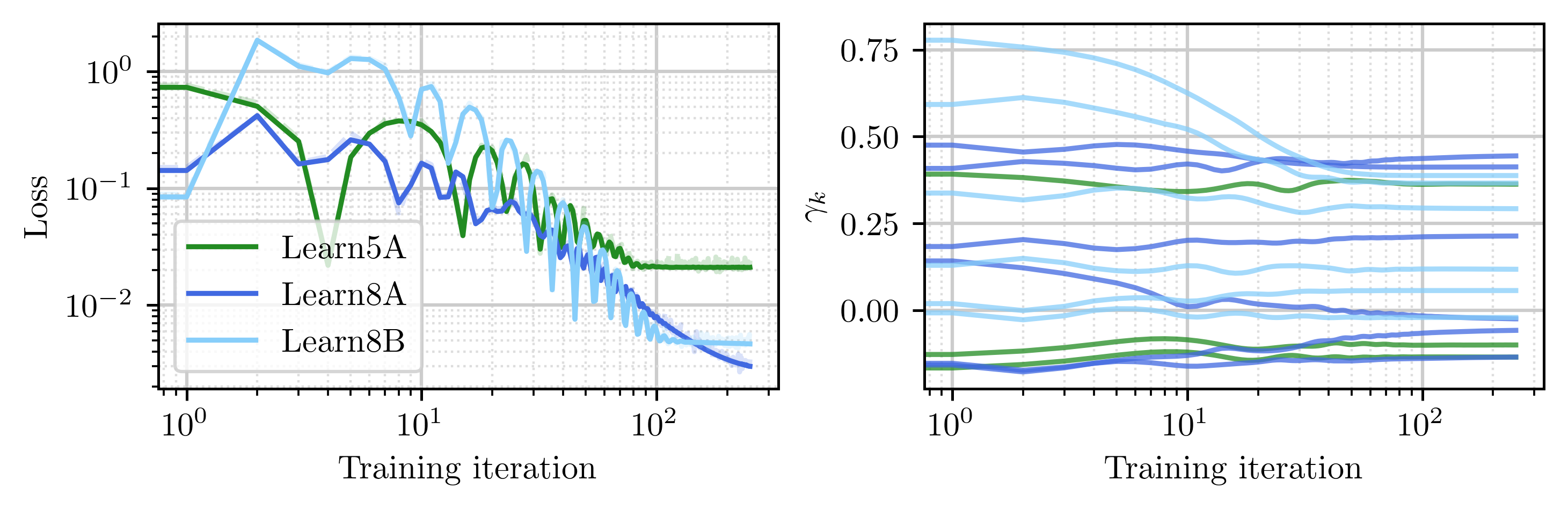}
    \caption{Training- and validation loss (left) and evolution of the splitting parameters for the three learned methods Learn5A, Learn8A and Learn8B where we used Adam for the stochastic optimisation.
    \label{fig:paramOptim}}
\end{figure}

Figure~\ref{fig:paramOptim} (left) shows the evolution of the training- and validation loss function in \eqref{eq:transLossFnFinite} for the Schrödinger equation \eqref{eq:potKinSchroedingerPDE} using  Adam for the stochatic optimisation in line 8 of the learning pipeline in Algorithm~\ref{alg:pipeline}. The loss is plotted as a function of the  iterations  for learning the optimal splittings with lengths $\splitDisc=5$ and $\splitDisc=8$. For $\splitDisc=8$, two different minima were identified, and these are labelled ``Learn8A'' and ``Learn8B'' respectively; the learned splitting method for $\splitDisc=5$ is referred to as ``Learn5A''. Both the training loss and the validation loss (which is much smoother since it is always evaluated on the same batch) are shown. The corresponding evolution of the splitting parameters $\allParams$ is shown in Figure~\ref{fig:paramOptim} (right). As the plot demonstrates, all parameters converge after 250 iterations.
\subsection{Comparison of different optimisers}\label{sec:optimiser_comparison}
Many different optimisers such as Adam \cite{kingma2014a}, AdaGrad \cite{duchi2011adaptive} or Lion \cite{chen2024symbolic} could be chosen in the SO step in line~8 of Algorithm~\ref{alg:pipeline}. Since the number of learned parameters is low, this includes second-order methods such as Levenberg-Marquardt \cite{levenberg1944method,marquardt1963algorithm} which require the evaluation of the Hessian. For each of these methods performance can be tuned by varying the relevant hyperparameters such as the learning rate, momentum and weight decay in the stochastic gradient descent on batched training data. 

\begin{figure}[tb]
    \centering
    \includegraphics[width=0.85\textwidth]{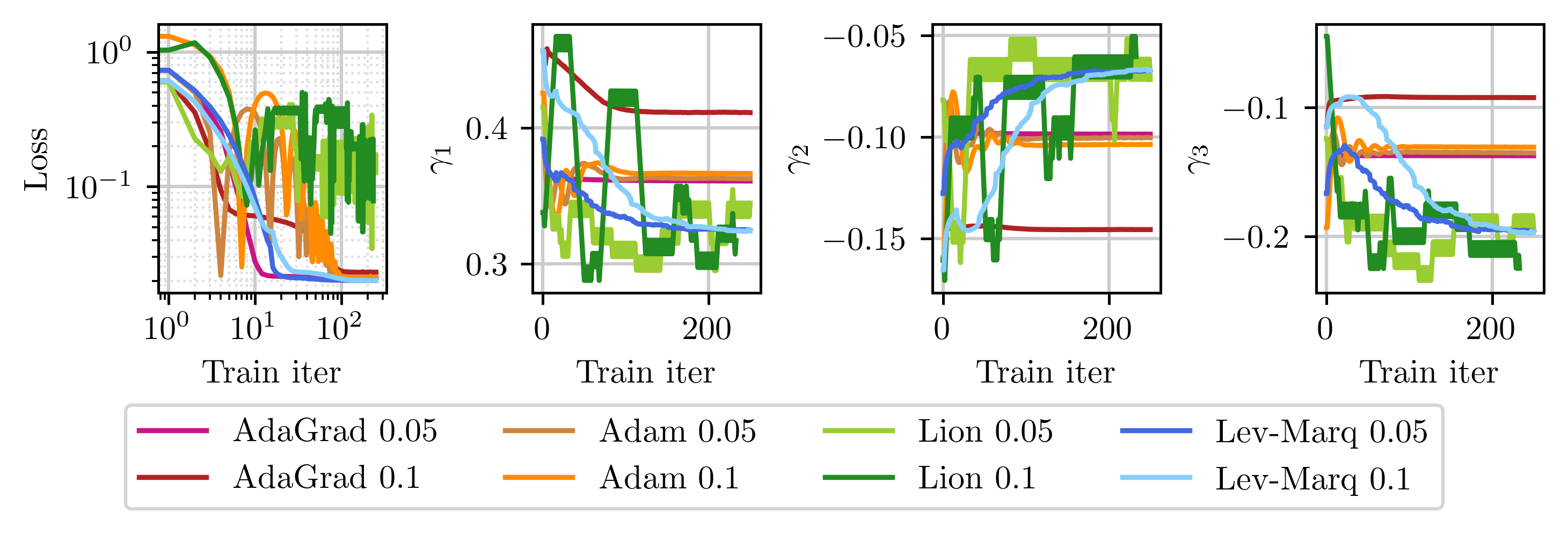}
    \caption{Comparison of different optimisers with associated learning rates 0.05 and 0.1 for finding a splittting method with three free parameters. The plot in the top left corner shows the evolution of the loss function; the other three plots visualise the corresponding evolution of the parameters $\gamma_1$, $\gamma_2$ and $\gamma_3$.
    \label{fig:allOptims}}
\end{figure}

Figure~\ref{fig:allOptims} illustrates the performance of AdaGrad, Adam, Lion and Levenberg-Marquardt on our training data when trying to find the minimum for a splitting method with three free parameters $\allParams_1$, $\allParams_2$ and $\allParams_3$.
To obtain this figure, we started near the novel learned splitting Learned5A by setting the initial values of the parameters to $\allParams_1=0.374129\pm \delta$, $\allParams_2=-0.109994\pm \delta$, $\allParams_3=-0.123223\pm \delta$ where $\delta\in \{-0.05,+0.05,-0.1,+0.1\}$ and explored the loss landscape. The loss function in the upper left plot in Figure~\ref{fig:allOptims} shows that all optimisers apart from Lion, converged to numerical methods with error of around $0.02$. For reference, the exact final values of the parameters obtained with the different methods and the corresponding value of the loss function are tabulated in Table~\ref{tab:final_parameter_values}. We conclude that they are all good candidates for optimisers to be used in line 8 of Algorithm~\ref{alg:pipeline}. As Adam is known to be an efficient and widely used algorithm, all other numerical experiments employ Adam.
\begin{table}[htb]
    \centering
    \begin{tabular}{|c|c|c|c|c|c|}
        \hline
        \textbf{optimizer} & LR & $\gamma_1$ & $\gamma_2$ & $\gamma_3$ & loss $\mathcal{L}$ \\
        \hline\hline
        \multirow{2}{*}{AdaGrad}  & 0.05 & 0.3609 & -0.09857 & -0.1375 & 0.02096 \\
        & 0.1 & 0.4114 & -0.1455 & -0.09246 & 0.2316 \\\hline
        \multirow{2}{*}{Adam} & 0.05 & 0.3627 & -0.1003 & -0.1353 & 0.02106 \\
        & 0.1 & 0.3665 & -0.1036 & -0.1309 & 0.0212 \\\hline
        \multirow{2}{*}{Lion} & 0.05 & 0.3446 & -0.07162 & -0.1938 & 0.1275 \\
        & 0.1 & 0.3169 & -0.06039 & -0.2246 & 0.07843 \\\hline
        Levenberg- & 0.05 & 0.3245 & -0.0671 & -0.1963 & 0.02011 \\
        Marquardt & 0.1 & 0.3242 & -0.06704 & -0.1966 & 0.2013 \\
        \hline
    \end{tabular}
    \caption{Final values of the parameters $\gamma_1$, $\gamma_2$ and $\gamma_3$ and loss function when trained over 250 epochs with different optimisers and learning rates LR.}
    \label{tab:final_parameter_values}
\end{table}
\section{Visualisation of splitting methods}\label{sec:splitting_vsualisation}
In the following we describe how the visualisation of splitting methods described in Section~\ref{sec:learning_longer_splittings} arises from the solution of a the simple two-dimensional IVP $\dot{\stateVar}(\timeVar) = [1,1]^\top$, $\trueState{0} = [0,0]^\top$, $\timeVar \in [0,1]$ for $\stateVar(\timeVar)=(\stateVar_1(\timeVar),\stateVar_2(t))^\top\in\mathbb{R}^2$.
The exact analytical solution $u(t)=[t,t]$ corresponds to a diagonal line from $[0,0]$ to $[1,1]$. The natural splitting is $\difEqComp{1}=[1,0]^\top$, $\difEqComp{2}=[0,1]^\top$. The corresponding sub-flows $\exactFlowComp{1}{\delta\timeVar}:(\stateVar_1(\timeVar),\stateVar_2(\timeVar))^\top\mapsto (\stateVar_1(\timeVar+\delta\timeVar),\stateVar_2(\timeVar))^\top$ and $\exactFlowComp{2}{\delta\timeVar}:(\stateVar_1(\timeVar),\stateVar_2(\timeVar))^\top\mapsto (\stateVar_1(\timeVar),\stateVar_2(\timeVar+\delta\timeVar))^\top$ translate the state parallel to the coordinate axes. For $0\le\tau \le 1$ we can now define the flux
\begin{equation}
    \numFlowCompOneStep{\difEq}{\tau\stepSize}(\, \cdot \, ; \potParams, \kinParams) = 
    \begin{cases}
        \exactFlowComp{1}{(\tau-S_{\splitInd'})\stepSize}\bigcirc_{\splitInd = 1}^{\splitInd'-1} \exactFlowComp{2}{\kinParams_\splitInd \stepSize} \circ \exactFlowComp{1}{\potParams_\splitInd \stepSize}&\text{for $S_{\splitInd'}\le\tau\le S_{\splitInd'}+\potParams_{\splitInd'}$}\\[1ex]
        \exactFlowComp{2}{(\tau-\potParams_\splitInd'-S_{\splitInd'})\stepSize} \circ \exactFlowComp{1}{\potParams_\splitInd' \stepSize} \bigcirc_{\splitInd = 1}^{\splitInd'-1} \exactFlowComp{2}{\kinParams_\splitInd \stepSize} \circ \exactFlowComp{1}{\potParams_\splitInd \stepSize}&\text{for $S_{\splitInd'}+\potParams_{\splitInd'}\le\tau\le S_{\splitInd'+1}$}
    \end{cases}
    \label{eqn:tau_flux}
\end{equation}
where the partial sums $S_{\splitInd'}$ are given by
\begin{equation*}
    S_{\splitInd'} = \sum_{\splitInd=1}^{\splitInd'-1} \kinParams_\splitInd + \potParams_\splitInd \qquad \text{for $\splitInd' = 1, 2, \dots, \splitDisc$ }.
\end{equation*}
As $\tau$ increases from $0$ to $1$, this flux will translate a given state $\stateVar$ to $\numFlowCompOneStep{\difEq}{\stepSize}(\stateVar; \potParams, \kinParams)$ along a connected path of straight line segments that are parallel to the coordinate axes and whose (orientated) lengths are given by $\kinParams_\splitInd$ and $\potParams_\splitInd$ respectively. This path can thus be regarded as a visualisation of the numerical time-stepping method defined by $(\potParams,\kinParams)$ which resolves the evolution within a single timestep of size $\stepSize$.
Figure \ref{fig:splitVisLearn} visualises the evolution of the initial state $[0,0]$ under different  learned and reference methods for $\stepSize=1$.
Consistency of a method is equivalent to the path ending at $[1,1]$, and for consistent methods symmetry is akin to having discrete $180^\circ$ rotational symmetry about the point $[0.5,0.5]$. By comparing Strang with $4 \times$ Strang in Figure~\ref{fig:splitVisLearn} it can be seen that decreasing the step size of a method brings the path closer to the diagonal line which represents the analytic solution at the cost of increasing the number of line segments.

\end{document}